\documentclass[10pt]{article}
\usepackage{tikz-cd}

\usepackage{float}
\usepackage{graphicx}  
\usepackage{epsfig}
\usepackage{amsfonts}
\usepackage{amssymb,amsmath,multicol}
\newtheorem{remark}{\HL{Remark}}
    \newcommand{\REMARK}[1]{ \begin{remark} #1 \end{remark}} 
	\usepackage{xcolor}

 \newcommand{\HL}[2][red]{#2}     %% CLEAN
		\newcommand{\COMMENT}[1]{}
\newcommand{\R}{ \mathbb R}  %% vectors
	\newcommand{\V}[1]{ {\boldsymbol{#1}}}  %% vectors
	\newcommand{\Vdot}[1]{ {\boldsymbol{\dot #1}}}

	\newcommand{\Vtilde}[1]{ {\boldsymbol{\tilde {#1}}}}
	\newcommand{\MAT}[1] {\begin{bmatrix} #1 \end{bmatrix}}
	\newcommand{\ARRAY}[1] {\begin{matrix} #1 \end{matrix}}
\usepackage{cancel}

\title{A Basic \HL{Mechanical and} Geometric Framework for \\ Quasi-Static Manipulation}

\author{Domenico Campolo$^1$, Franco Cardin$^2$ \\
	\small{ $^1$School of Mechanical and Aerospace Engineering,} \\
	\small{Nanyang Technological University (NTU), Singapore.} \\
	\vspace{3mm} 	\small {\tt d.campolo@ntu.edu.sg} \\ 
	\small{ $^2$Dipartimento di Matematica “Tullio Levi-Civita”,} \\ 
	\small{ 	Universit\`{a} degli Studi di Padova,  Italy}\\
	\small {\tt cardin@math.unipd.it} \\
}
\date{}
\begin{document}
\maketitle

%\tableofcontents
\abstract 

In this work, we propose a geometric framework for analyzing 
mechanical manipulation, for \HL{instance}, by a robotic agent. Under the assumption of conservative forces and quasi-static manipulation, we use energy methods to derive a metric. 

\HL{In the first part of the paper, we}
review how \HL{quasi-static mechanical manipulation tasks can be naturally described via the so-called \textit{force-space}, i.e. }
the cotangent bundle \HL{of the configuration space,} and its Lagrangian submanifolds. 
Then, via a second order analysis, we derive the \HL{control} Hessian of total energy. As this is not necessarily positive-definite, from \HL{an optimal} control perspective, we propose the use of the squared-Hessian, also motivated by insights derived from both mechanics (\HL{Gauss'} Principle) and biology (Separation Principle).

\HL{In the second part of the paper, we apply } such methods to \HL{ the problem of an elastically-driven, inverted pendulum. }
\HL{
Despite its apparent simplicity, this example is representative of an important class of robotic manipulation problems for which we show how a smooth elastic potential can be derived by regularizing mechanical contact. 
%Via standard numerical methods over a grid, we derive a discrete set of equilibria. 
We then show how graph theory can be used to connect each numerical solution to `nearby' ones, with weights derived from the very metric introduced in the first part of the paper. 
%In summary, after numerical discretization, the original manipulation problem reduces to finding  a minimum paths over a weighted graph.
}

~

\newpage

\section{Introduction}

The definition of \textit{manipulation}, as found in dictionaries, invariably contains two elements: the first, relating to the act of operating either by hand or by mechanical means on the environment; the second, relating to skillfulness or to the purpose of gaining some advantage. Therefore, in its basic definition, manipulation involves not only \textit{physical interaction }of an agent with the surrounding environment but also a sense of \textit{quality} or \textit{cost/benefit}. 
\HL
{
Along these lines, this work proposes an \textit{optimality measure} for mechanical manipulation, with focus on \textit{quasi-static manipulation} of network of elastically connected rigid bodies, with a finite number of degrees of freedom (dof).

\subsection{Geometric setting: the force-space}
The starting point for the analysis of a mechanical system is the definition of its \textit{configuration space}, conventionally defined as the space in which motion takes place, typically involving different sets coordinates and identified with some differentiable manifold $Q$, e.g. see \cite{Abraham and Marsden 2008, Arnold 1989, Bullo and Lewis 2005}. 
In this work we shall only consider systems characterized by a finite number $n$ of degrees of freedom,  leading to an $n$-dimensional configuration space, i.e. $\dim(Q)=n$.  
Many systems of interest will fall in this category,  including robots and  environment they interact with \cite{Mason 2018}.

To describe the evolution of a system over its configuration space, classical mechanics typically presents a dual perspective: $i)$ the \textit{Lagrangian} formulation on the \textit{tangent bundle} $TQ$, the space of positions and velocities, leading to the Euler-Lagrange (EL) equations; $ii)$ the \textit{Hamiltonian} formulation on the \textit{cotangent bundle} $T^*Q$, leading to Hamilton's equations, e.g. see \cite{Cardin 2015, Benenti 2011}.
Both tangent and cotangent bundles are  $2n$ dimensional manifolds, i.e. $\dim(TQ)=\dim(T^*Q)=2n$. 

Traditionally, the mechanical interpretation of the cotangent bundle is that of \textit{phase-space} where, at each configuration $\V q\in Q$,  linear elements of the $n$-dimensional cotangent space $T_{\V q}^*Q$, or \textit{covectors}, are identified with the generalized mechanical \textit{impulse}. 
However, as also highlighted in Ref. \cite{Tulczyjew 1989}, \cite[\S 4.1]{Benenti 2011}, there is a second mechanical interpretation, whereby a covector  $\V f\in T_{\V q}^*Q$  can be viewed as a generalized mechanical \textit{force} which, via natural  pairing $\langle\cdot, \cdot \rangle : T_{\V q}^*Q\times T_{\V q}Q\to\R$ with displacements $\delta \V q\in T_{\V q}Q$, returns \textit{work} $\langle\V f, \V \delta\V q\rangle$.  This \textit{force-space} interpretation finds applications ranging  from control of static mechanical systems to thermodynamics \cite[Chapter 8]{Benenti 2011} and will be provide the basic geometric framework for this work. 

Although the Lagrangian formulation is predominant in the robotics literature  \cite{Bullo and Lewis 2005, Lynch and Park 2017}, we find that, as mechanical manipulation inherently involves interaction forces, the force-space interpretation of the cotangent bundle  provides a natural setting for the analysis of manipulation tasks. Furthermore, as a starting point for quasi-statics analysis, we shall be interested in the \textit{subset of equilibria}. To this end, consider a mechanical spring undergoing (one-dimensional) extension or compression. Its possibly nonlinear response is typically described by a force vs. displacement graph $\{(q, f(q))\in\R^*\times \R\simeq \R^2 : q\in \R \}$. 
In the simplest cases, this is given as an explicit analytical relationship $q\mapsto f(q)$ but it is not hard to find examples of simple configurations of linear springs leading to  multivalued maps, e.g. see \cite[\S8.1]{Benenti 2011}.  
The key point, here, is that although cotangent bundles represent a natural $2n$-dimensional ambient space for \textit{any} subset of force-displacement pairs, the actual description of a \textit{specific }system (e.g. a  spring, as in the example above) lies in definition of a \textit{specific }subset of the cotangent bundle. Under certain technical conditions and for the cases of interest in this work, these subsets are actually $n$-dimensional submanifolds of the cotangent bundle,  known as \textit{Lagrangian submanifolds} 
\cite{Tulczyjew 2011, Tulczyjew 1989, Tulczyjew 1982, Benenti 2011, Cardin 2015}. 

\subsection{Quasi-static assumptions and coordinates as controls}
For a large variety of manipulation tasks, including mechanical assembly by a robotic or a human agent, inertial forces are often negligible, suggesting a \textit{quasi-static analysis} of the problem \cite{Whitney 1982, Mason 2001, Mason 2018}. This type of approximation is also known as \textit{overdamped regime} and there is a close connection between classical mechanics, minimum principles and non-equilibrium thermodynamics. For a recent account, the reader is referred to \cite{Podio-Guidugli and Virga 2023}. When inertial forces become negligible, the Hamiltonian (as well as the Lagrangian, with a change of sign) simply accounts for the potential energy of the system $W:Q\to \R$. 
Dissipation can also be accounted for via, for example, \textit{Rayleigh's dissipation potential} $R(\V q, \Vdot q)$, a positive semidefinite quadratic form for all configurations $\V q$ \cite{Podio-Guidugli and Virga 2023}.

The specific type of evolution, which depends on the type of dissipation, is beyond the scope of this work. What is important in our discussion is that, under these  assumptions, $i)$ the system will evolve and eventually settle at some stable equilibria (as any gradient system); $ii)$ the presence and type (stable/unstable) of equilibria will only depend on the potential and not on the dissipation function. 

So far, only the geometric setting and the evolution of mechanical systems has been discussed but \textit{how can it  be controlled?}
As reviewed and discussed in \cite{Bressan and Rampazzo 2010}, there are two fundamentally different  ways to control mechanical systems: $i)$ by exerting additional forces; $ii)$ by directly controlling some the coordinates of the system. Although the former represents the most common approach, there is a rich and consolidated literature on the theory of controls operated by coordinates in a fully dynamic environment, from the pioneering works of Bressan\footnote{It should be noted that Ref.~\cite{AldoBressan88, AldoBressan1989} should be attributed to Aldo Bressan while Ref.~\cite{Bressan and Rampazzo 2010, Bressan and Motta 1993, Bressan and Rampazzo 2010, Bressan et al 2013} should be attributed to Alberto Bressan.} \cite{AldoBressan88, AldoBressan1989} and  Marle \cite{Marle 1990, Marle 1991} and with  a very active and still operating school \cite{Bressan and Motta 1993, Bressan and Rampazzo 2010, Bressan et al 2013, Cardin and Favretti 1998}.
Along this avenue, the approach adopted in this work is to consider \textit{some} of the configuration variables as directly controlled by an  external agent and, in particular, splitting the variables into \textit{internal states} ($\V z$) and \textit{control inputs} ($\V u$), i.e. $\V q\equiv (\V z, \V u)$. 

The problem is therefore that of finding suitable control inputs $t \mapsto \V u(t)$ capable of  driving the system from a given initial configuration $\V z_0$ to a desired final configuration $\V z_1$. Furthermore, we shall operate under \textit{quasi-static regime}, i.e. we shall assume  that the control inputs $\V u(t)$ are varied so slowly with respect to the system dynamics that, at any time $t$, we can expect the full state of system $(\V z(t), \V u(t))$ to be always at equilibrium. 
This is clearly an approximation and it should be noted that the variable $t$ is simply a parameter and that time plays no role in this quasi-static regime.

\subsection{Optimality}

It is well known that in robotics, just like in more general mechanical systems, the inertial term (i.e. the term accounting for the kinetic energy of the system) embeds a \textit{Riemannian metric} \cite{Bullo and Lewis 2005, Zefran and Bullo 2005}. 
However, in an overdamped regime, with negligible inertial effects, \textit{what can possibly play the role of a metric?}

Since the 80s, Loncaric \cite{Loncaric 1987, Loncaric 1991} hinted at \textit{stiffness} as a possible source of positive, definite matrices, in addition to inertia and damping. The analysis was carried out on the Special Euclidean group of rigid body motions $SE(3)$. Elastic interactions were considered via (scalar) potentials $W:SE(3)\to \mathbb R$ and their elastic \textit{stiffness} matrix, defined as its \textit{Hessian} at \textit{critical points} $dW=0$ (the only points where Hessians behave tensorially) \cite{Loncaric 1991}. It should also be noted that \textit{compliance/impedance control} is a prominent approach in modern robotics to control physical interaction, e.g. see \cite{Hogan 2022}.

However, when considering a controlled system with internal states $\V z$, control inputs $\V u$ and potential $W(\V z, \V u)$, while the positive-definiteness of the Hessian with respect to internal variables is linked to stability  the system itself, as shown later, the \textit{control Hessian} 
$\V G_m(\V u)$ in Eq.\eqref{EQ:Gm(u)}, i.e. the second-order mixed derivatives with respect to control variables at equilibrium points, needs not be.
Inspired by Gilmore \cite{Gilmore 1981},  the approach taken  in this work bears some similarity  with thermodynamics but also a major difference: due to the specific nature of thermodynamic variables, a consequence of the second principle of thermodynamics is that  the corresponding control Hessian in the thermodynamic setting is always positive-definite \cite[\S10-12]{Gilmore 1981}.
 
To overcome the non-positive definiteness of the control Hessian in our framework, one of the main contributions of this work is to show how the \textit{squared control Hessian}  can be used as a \textit{manipulation metric}.
It should be noted that the idea of a squared-Hessian as a \textit{Riemannian metric} is not new. In fact, Hermosilla \cite{Hermosilla 2016} shows how a metric could be derived from the squared-Hessian of a given optimization cost, however, without specifying what this cost might be. In our work, we show how the very potential $W(\V q, \V u)$ defining the system, naturally leads to a metric via the squared-Hessian.

While the squaring of the Hessian might appear as a choice of convenience, we also show that this is  motivated by the  \textit{Separation Principle}, a biological insight posing that the human brain processes static (or configuration-dependent) and dynamic (or velocity-dependent) force fields \textit{separately}, as suggested by experimental evidence \cite{Atkeson and Hollerbach (1985), Guigon et al (2000), Hollerbach and Flash (1982), Kurtzer et al (2005), Tommasino and Campolo (2017)}. In brief, as it will be detailed later, rather than optimizing the `total effort' (e.g. the muscular activation) to move a limb from an initial to a final location, the human brain appears to optimize the `net effort', i.e. the total effort minus the effort needed to maintain equilibrium. Imagine, for example, writing on a whiteboard: while most of the effort is needed to sustain the weight of the arm, much less force (net effort) is needed to actually move the arm to draw letters.

Although biologically inspired, the Separation Principle joins a class of local principles  such as Gauss' Principle \cite{Gauss 1829}, an ancestral tool for describing dynamics of  ideal, constrained mechanical system. Born in a genuine non-static environment, it states that along the admissible motions, point by point
in the host tangent space, \textit{the square of the norm of the global reaction forces has to be minimum for any admissible variations of the accelerations}. This notable Principle, described in an intrinsic geometrical framework \cite{Cardin and Zanotto 1989}, has been often utilized in many contexts: statistical mechanics, see \cite[p. 320]{Gavallotti 1999}, and fluid dynamics,  \cite[p. 444, p. 451]{Gavallotti 2013}; furthermore, by thinking of our aims, in various robotic settings, e.g. \cite{Lilov and Lorer 1982}.

\subsection{Numerical Approaches}

From the discussion above, it emerges that a geometric description of mechanical manipulation might be given by defining specific submanifolds (of an appropriate ambient space) and a metric on them. However, as it will be detailed later, such manifolds are often only \textit{implicitly }defined. Resorting to numerical analysis, the only available description would then be a set of $m$ \textit{samples} $\{(\V q_i, \V f_i)\in T^*Q\}_{i=1}^m$ (we shall not consider the problem of noisy information, here). The issue is that a discrete set of points, although representative, is not the implicit submanifold we are looking for. 
This problem is akin to\textit{ manifold learning}, which is a vast topic and we shall refer the reader to a recent account \cite{Joharinad and Jost 2023}. A standard  approach to manifold learning is based on \textit{Graph Theory}, whereby sampled data in a high-dimensional ambient space represent \textit{nodes} of the graph and the connection between edges determines the topology of the graph.

Unlike typical manifold learning problems, where an  underlying manifold is assumed but not even its dimensionality is known, let alone a metric, in our case we assume that a careful analysis of the problem (starting from the definition of an appropriate configuration space) is available. While the submanifolds of interest are only defined implicitly, once a sample is given (e.g. numerically determined), geometric constructions such as tangent space and quadratic form of the metric, at the sample itself,
can be explicitly (albeit numerically) determined. As detailed below, we can then construct the connectivity of a graph by determining which samples are neighbors, associating a weight to the connection between two neighboring samples via the squared-Hessian metric.

\subsection{Organization of the paper}
This paper mainly consists of two parts. In the first part, after recalling  basic definitions for Lagrangian submanifolds, Sec.~\ref{sec:EM} will introduce the Equilibrium Manifold, in fact a submanifold of the configuration space and its first order approximation. A second order approximation will be introduced in Sec.~\ref{sec:metric} and this will be used to derive a metric from the squared-Hessian. In this section a clear connection between Separation Principle and squared-Hessian will be outlined.

In the second part of the paper, Sec.~\ref{sec:toy_model}, the theory developed in the first part will be applied to a `toy model': an inverted pendulum elastically driven by an agent.
On the one hand, stabilization of an inverted pendulum is a classical benchmark in control theories, with one of the earliest accounts dating back to Kapitza \cite{Kapitza 1951}.
On the other hand, despite its apparent simplicity, this example is representative of an important class of robotic manipulation problems. In particular, we show how a smooth elastic potential can be derived by regularizing mechanical contact. This is essential in order to derive a differentiable model.
Finally, given a smooth but highly nonlinear potential, standard numerical methods can be used to derive a discrete set of equilibria. Graph theoretical methods will then be deployed to construct a weighted graph which connects `nearby' points based on the smooth theory developed in the first part of the paper, in particular the squared-Hessian metric. 
The original manipulation problem is then reduced to finding a minimum paths over a weighted graph.

\subsection{Notation}
In terms of notation, both vectors and covectors will be represented as column arrays and denoted with bold, lower-case symbols e.g. $\delta \V q=[\delta q^1, \dots, \delta q^n]^T\in\R^n$ and  $\V f=[f_1, \dots,  f_n]^T\in(\R^n)^*$, where $(\cdot)^T$ denotes matrix transposition. 
The pairing $\langle\V f, \delta \V q\rangle$ can therefore be evaluated as the scalar  $\V f^T\delta\V q$.
%Bold, upper-case symbols will be used to denote matrices, e.g. $\V A\in\R^{n\times m}$.

In local coordinates $\V q\in\R^{n}$, given a scalar potential potential $\V q\mapsto W(\V q)$, the components of its differential $dW$ in local coordinates will be denoted as $\nabla_{\V q}W\equiv [\partial_{q^1}W, \dots, \partial_{q^n}W]^T$, where the \textit{nabla} (column) operator is defined as  $\nabla_{\V q}:=[\partial_{q^1}, \dots, \partial_{q^n}]^T$.
We shall further introduce the shorthand notation $\nabla_{\V z\V z}^2\equiv \nabla^T_{\V z}\nabla_{\V z}$ and $\nabla_{\V u\V u}^2\equiv \nabla^T_{\V u}\nabla_{\V u}$, where  $\V z\in\R^N$ and $\V u\in R^K$,  for the Hessians as well as for the mixed-derivatives operators $\nabla_{\V u\V z}^2\equiv \nabla^T_{\V u}\nabla_{\V z}$ and $\nabla_{\V z\V u}^2\equiv \nabla^T_{\V z}\nabla_{\V u} =(\nabla_{\V u\V z}^2)^T$.

}

\section{Equilibrium Manifold - a first order analysis}
\label{sec:EM}
This section describes the first-order geometry of the set of equilibria for \HL{finite-dimensional} and conservative mechanical systems controlled by an external agent and operating under quasi-static assumptions. In particular, we will consider systems consisting of a finite number of rigid bodies and particles, interconnected via springs and possibly subjected to gravity. Some of the degrees of freedom will be considered as parameters directly controllable by an agent. 

In what follows, first, we shall provide basic definitions which, albeit standard in geometric mechanics, are here specialized to the static analysis. It will be apparent that the natural space for our analysis is \HL{the} cotangent bundle \HL{of the configuration space} (here referred to as force-space, in the context of statics) but not all of it, only its Lagrangian submanifolds.
Then, these general definitions will be specialized to the problem at hand, which involves controlled mechanical systems. This will correspond to a splitting of degrees of freedom of the configuration space into those which are directly controllable by an agent and those which are not.

\subsection{Force-space and its Lagrangian submanifolds}
Consider a \textit{conservative, \HL{finite-dimensional} mechanical system} characterized by a finite set of configuration variables $\V q\equiv \HL{[q^1, \dots, q^n]^T}$, at least within some open set  $\HL{Q}\subset \mathbb R^n$ as we shall mainly be concerned with \textit{local properties}.
Due to conservativity, assume the existence of a  $C^\infty$ scalar potential energy $W:\HL{Q}\to\mathbb R$. At any configuration $\V q\in \HL{Q}$, 
\HL{
the force arising from the given potential corresponds to (the negative of) its  \textit{differential} $-dW$, an element of the \textit{cotangent space} $T_\V{q}^*\HL{Q}$, and will be expressed in local coordinates as:
\begin{equation}\label{EQ:f=gradW}\nonumber
	\V f\HL{\equiv [f_1, \dots, \HL{f_n}]^T}=-\nabla_{\V q} W %\equiv [\partial_{q^1}W, \dots , \partial_{q^n}W]^T. %\in T^*_{\V q} Q
\end{equation}
}
From the onset, we shall highlight the {duality} between forces as {covectors}, i.e. elements of the cotangent space $T^*_{\V q}\HL {Q} $ at a given configuration $\V q$, and infinitesimal displacements $\delta \V q=\HL{[\delta q^1, \dots, \delta q^n]^T}$ as vectors, i.e. elements of the tangent space $T_{\V q}\HL{Q}$. This duality is reflected as  natural pairing in the definition of \textit{work:}

\begin{equation}\label{EQ:virtual_work}
	\delta \mathcal W:=\sum_{i=1}^nf_i\delta q^i\equiv \V f^T\delta \V q.
\end{equation}

Just like force vs. displacement graphs are often used to geometrically represent elastic properties of simple springs as curves in a plane, for a given potential $W$ one can consider the \textit{graph}:
\begin{equation}\label{EQ:Lambda^W}
	\Lambda^W = \{ (\V q, \V f)~:~\V q\in \HL{Q}, \ \ \V f= - \nabla _{\V q} W \}
	  \subset T^*\HL{Q} ,
\end{equation}
which can be shown to be an $n$-dimensional submanifold of the cotangent bundle $T^*\HL{Q}$, a $2n$-dimensional space hereafter referred to as \textit{force-space}.

\HL{Motivated by the definition \eqref{EQ:virtual_work},} one can introduce the so-called \textit{Liouville 1-form}
\begin{equation}\nonumber
	\Theta :=\sum_{i=1}^nf_i dq^i  , 
\end{equation}
a natural differential form on any cotangent bundle \cite[\S37-B]{Arnold 1989}.  
\HL{
In fact,  the most general 1-form on $T^*\HL{Q}$ would be $\sum_i a_i(\V q, \V f)dq^i+\sum_i b_i(\V q, \V f)df_i$ but, motivated by definition  \eqref{EQ:virtual_work}, one obtains $\Theta$ by setting \HL{$a_i=f_i$} and $b_i=0$, for $i=1, \dots, n$.
}

One can immediately verify that, when restricted to $\Lambda^W$, the Liouville 1-form reduces to an \textit{exact} differential  $\Theta|_{\Lambda^W} =dW$, which can be seen as an integrability condition.  
\HL{
A consequence of $\oint dW=0$ is that any loop $\V\gamma:[0 \ 1]\to \HL{Q}$, with $\V\gamma(0)=\V\gamma(1)$ on the configuration space, can be lifted on $\Lambda^W$, i.e. $\Vtilde\gamma:t \mapsto (\V\gamma(t), \nabla_{\V q}W|_{\V \gamma(t)})$ and
	$\oint_{\Vtilde\gamma}\Theta=0$. 
}
In other words, $\Lambda^W$ is the image of the differential $dW$, i.e. $\Lambda^W=\HL{\text{im}(-dW)}$ and is described by its graph \eqref{EQ:Lambda^W}.
More in general, the graph of any \textit{closed} 1-form (not just exact differentials) defines a Lagrangian submanifold $\Lambda \subset T^*\HL{Q}$.

\HL{
	More formally, an $n$-dimensional submanifold $\Lambda\subset T^*\HL{Q}$ is a \textit{Lagrangian submanifold} if the following integrability condition holds 	
	 $$d\Theta|_{\Lambda}=0$$ 	where $d$ represents the exterior derivative and $d\Theta$ is the \textit{symplectic  2-form} naturally defined on any cotangent bundle \cite[\S37-B]{Arnold 1989}. 
	The equality $d\Theta|_{\Lambda}=0$ is equivalent to stating that `$\Theta_\Lambda$ is \textit{closed}' and is nothing but a `curl-free' condition. 
}

This leads one to define more general \textit{multi-valued} Lagrangian submanifolds $\Lambda$ \cite{Ekeland 1977}
generated by 
 so-called \textit{Morse families} $W$ \cite{Cardin 2015}, locally expressed as 

\begin{equation}\nonumber
\Lambda = \{(\V q, \V f) ~:~  \V f = - \nabla_{\V q}W(\V q, \V \xi), \ \V 0 = \nabla_{\V \xi} W(\V q, \V \xi)\},
\end{equation}
where now \HL{$W:\HL{Q}\times \mathbb R^K\to \mathbb R,(\V q, \V \xi)\mapsto W(\V q, \V \xi)$} is {thought of as a sort of} potential parameterized by an auxiliary $K$-dimensional variable $\V \xi \in \mathbb R^K$,
\HL{introduced by Maslov-Hormander  in their classical theorem characterizing the Lagrangian submanifolds \cite[\S2.3]{Cardin 2015}.}
For this parameterization to be possible, the following \textit{maximum-rank condition }is necessary \cite{Cardin 2015}

\begin{equation}\label{EQ:general_max_rank}
	\text{rank}\MAT{\nabla^2_{\V q\V \xi}W \ \nabla^2_{\V \xi\V \xi}W}_{\nabla_{\V \xi}W=\HL{\V 0}}= K:  max.
\end{equation}
\HL{a condition which allows us to apply the implicit function theorem and guarantees that $\Lambda$ is a smooth manifold.} This extra parameterization will play an important role in the analysis that follows.
 
\HL{Next, the geometric setting of force-space will be used to} describe a network of rigid bodies and particles, possibly interconnected via generalized springs and under the possible influence of gravity. 

\subsection{Interconnected mechanical systems}
Following \cite{Bullo and Lewis 2005}, the configuration space $Q$ of interconnected mechanical systems can be thought of as a $C^\infty$-immersed submanifold of the product space
\begin{equation}\label{EQ:Q_subset}\nonumber
	Q\subset \underbrace{SE(3)\times \dots \times SE(3)}_{N_b \text{ rigid bodies}}\times \underbrace{\mathbb R^3 \times \dots\times \mathbb R^3}_{N_p \text{ particles}}	
\end{equation}
of dimension $\dim Q \le 6 N_b+3 N_p$. We shall assume that $K\ge 1$ degrees of freedom (dof) are directly controllable by an agent while the remaining $N=\dim Q - K$ dof can only be indirectly influenced via generalized elastic forces or via external forces such as gravity. 
%In robotics, such  systems are often referred to as `\textit{underactuated}', meaning that the directly controllable dof ($K$) is strictly less than the total number of dof ($\dim Q$). 

%\subsubsection*{Local analysis}
The analysis that follows is mostly local, relying on the Implicit Function theorem and on Taylor's expansion as main analytical tools. 
In particular, we shall assume the existence of equilibria in some open neighborhood of $Q$ parameterized
	via a set of coordinates
$(\V z, \V u)\in \mathcal Z \times \mathcal U$, where $\mathcal Z\subset\mathbb R^N$ and $\mathcal U\subset\mathbb R^K$ are open sets. These coordinates represent

\begin{itemize}
	\item \textit{internal states }$\boldsymbol  z\in\mathcal Z\subset \mathbb R^N$, non directly controllable; %(i.e. under-actuated);
	\item  \textit{control inputs }$\boldsymbol  u\in\mathcal U\subset \mathbb R^K$, directly controllable by the agent, therefore assumed as an input to the system.
\end{itemize}

We shall further assume that the system will solely be subjected to conservative forces (either via gravity or internal springs) and, in particular, the existence of a smooth ($C^\infty$) potential energy 
\begin{equation}\label{EQ:W}\nonumber
	W(\boldsymbol z, \boldsymbol  u):\mathcal Z\times \mathcal U\to \R,
\end{equation}
which will play the role of \textit{Morse family}, as defined above and  \HL{variables renamed as}  
$(\V q, \V \xi) \to (\V u, \V z).$    
  
\HL{
\REMARK{
In the change of variables $(\V q, \V \xi) \to (\V u, \V z)$, the role of $\V q$ will be played by the control inputs  $\V u$ while the role of auxiliary variables $\V \xi$ will be played by the remaining variables $\V z$, i.e. the internal states.
This appears to be an important change of perspective with respect to standard interpretation of these objects.
}
}
Considering the control inputs $\V u$ as parameters, mechanical equilibria correspond to stationary points of the potential with respect to internal variables $\V z$, which is equivalent to a zero-condition for \textit{internal forces} $\nabla_{\V z} W$.  \HL{With perhaps an abuse of notation, the wording ``internal forces'' really refers to differentials with respect to ``internal'' variables $\V z$, this also includes gravity which is definitely not internal to the system.}
\HL{Also, as it will be detailed in Sec.~\ref{sec:Equilibrium_Forces}, while at equilibrium internal forces are expected to null, i.e. $\nabla_{\V z}W=\V 0$, in general control forces are now, i.e $\nabla_{\V u}W\neq \V 0$.} 

For a given control input $\boldsymbol u^*$, there exist possibly multiple solutions $\boldsymbol z_m^*$ (with $m\ge 1$ being an integer denoting the \textit{multiplicity} of  equilibria) to the \textit{equilibrium equation}:

\begin{equation}\label{EQ:equil}
	\nabla_{\V z}W (\V z_m^*, \V u^*)=\V 0,
\end{equation}
where $\nabla_{\V z}$ denotes the column operator  $[\partial_{z^1}\ \dots \ \partial_{z^N}]^T$ and, similarly, $\nabla_{\V u}\equiv [\partial_{u^1}\ \dots \ \partial_{u^K}]^T$. 

Once the input $\V u^*$ is fixed, stability of the mechanical system is determined by the positive-definiteness of the Hessian $\nabla_{\V z\V z}^2W$ \cite{Bazant and Cedolin 2010}. If the Hessian is full-ranked at a given equilibrium, i.e.
\begin{equation}\label{EQ:max_rank_Wzz}
	\hspace{1cm}	\text{rank}(\nabla^2_{\V z\V z}W|_{\HL{(\V z_m^*, \V u^*)}}) = N \hspace{1cm} (\textit{max-rank)},
\end{equation}
then such an equilibrium will be referred to as {\em non-critical}.

Assuming the existence of {(at least)} one solution $(\V z_m^*, \V u^*)$ and max-rank condition \eqref{EQ:max_rank_Wzz}, then by the Implicit Function Theorem, here restated as in \cite[Th. 2-12]{Spivak 2018}, there exist an open $\mathcal U'\subset \mathcal U$ containing $\V u^*$, an open set $\mathcal Z'\subset \mathcal Z$ containing $\V z^*$ and a differentiable map $\V z_m :\mathcal U' \to \mathcal Z'$ such that $\V z_m(\V u)\in \mathcal Z'$ is \textit{unique} and $\nabla_{\V z}W(\V z_m(\V u), \V u)= \V 0$, for all $\V u\in \mathcal U'$.

\HL{
\subsection{Local branches of Lagrangian submanifolds}
}
Locally, in an open neighborhood $\mathcal Z'\times \mathcal U'$ of a non-critical equilibrium $\V z_m^*, \V u^*$, the \textit{Equilibrium Manifold} (EM) can be described as a $K$-dimensional submanifold
	\begin{equation}\label{EQ:EM}
		\text{EM}=\{(\V z, \V u)\in{\mathbb R^N\times \mathbb R^K}: \nabla_{\boldsymbol z}W{(\V z, \V u)}=\V 0\}.
	\end{equation}
For each $\V u\in\mathbb R^K$, the values $\V z\in\mathbb R^N$ are typically  \textit{discrete } and EM 
can be parameterized by controls $\V u$, at least locally,
in an open neighborhood of a non-critical equilibrium 
$(\V z_m^*,\V u^*)$.
{EM is simply a geometric feature of the Lagrangian submanifold $\Lambda$ generated by $W(\V z, \V u)$.}

In summary, for a given control input $\V u^*$, the equilibrium equation \eqref{EQ:equil} can be solved (e.g. numerically) to derive possibly multiple solutions $(\V z_m^*, \V u^*)$, with \HL{multiplicity} $m=1,2, \dots$. Solutions $(\V z_m^*, \V u^*)$ satisfying  the max-rank condition \eqref{EQ:max_rank_Wzz} 
\HL{also satisfy the condition $\det(\nabla^2_{\V z\V z}W|_{(\V z_m^*, \V u^*)})\neq 0$. Therefore, around each such solution,}
there exists an open neighborhood $\mathcal Z'\times \mathcal U'$ on which the EM is a submanifold and we shall refer to it as \textit{local branch} of EM.
It can be shown that non-critical equilibria are \textit{isolated} and can, at least in principle, be controlled (sliding along the EM) by varying the control parameters $\V u$. 

To this end, let us consider curves on the control space $s \mapsto \V u(s)\in \mathcal U$ and their lift
on $\mathcal Z\times \mathcal U$ 
\begin{equation}\nonumber
	s \mapsto (\V z_m(\V u(s)), \V u(s)) \in \text{EM}
\end{equation}
on the equilibrium manifold, parameterized by a scalar $s\in \mathbb R$, at least locally, on open intervals including $s=0$, for which $\V u(0)=\V u^*$ and $\V z_m(\V u(0))=\V z_m^*$
%where the asterisk $^*$ reminds us that $(\V z_m^*(\V u^*(s)), \V u^*(s))$ 
For every such a curve, the equilibrium condition \eqref{EQ:equil} is satisfied for all $s$, i.e. 
\begin{equation}\label{EQ:equi_lambda}
	\nabla_{\V z}W(\V z_m(\V u(s)), \V u(s))=\V 0.
\end{equation}
Differentiating \eqref{EQ:equi_lambda} with respect to $s$, leads to
\begin{equation}\label{EQ:equil_delta}
	\nabla_{\V z\V z}^2W \delta\, \V z + \nabla^2_{\V u\V z}W \,\delta \V u=\V 0,
\end{equation}
\HL{where $\delta \V u:= \dot{\V u}\,ds$ and 
$\delta \V z:=\nabla^T_{\V u}\V z_m \delta\V u$ 
have been introduced to maintain similarity with the nomenclature as in \cite[Sec. 5]{Gilmore 1981}. 
}
Thanks to the max-rank condition of $\nabla^2_{\V z\V z}W$,  \eqref{EQ:equil_delta} can be rewritten as 
\begin{equation}\label{EQ:1st_order}
	{\delta \V z = - (\nabla_{\V z\V z}^2W)^{-1}\nabla^2_{\V u\V z}W\delta \V u. }
\end{equation}

\REMARK{
    The tangent vector variation $\delta \V z$ represents the response to a change of   macroscopic internal equilibrium $\V z$ under the modification $\V u$ of the controls; all this, on the branch marked by $\V z_m$.
    }

We are now able to express a first order approximation $\Pi_{*_m}$ to the equilibrium manifold \eqref{EQ:EM}, i.e. $T_{*_m}\text{EM}$ at a given equilibrium point $*_m\equiv (\V z_m^*, \V u^*)$ as

\begin{equation}\label{EQ:EM_explicit}
	\Pi_{*_m} 
	= \{(\V z_m^*  \underbrace{- (\nabla^2_{\V z\V z}W)^{-1}\nabla^2_{\V u\V z}W\,\delta \V u}_{\delta \V z}, \V u^*+\delta \V u )
	~: ~ \delta \V u\in \mathbb R^K\}	.
\end{equation}	

As a final note, in this work we shall always assume stability of the mechanical system as well positive definiteness of $\nabla^2_{\V z \V z}W$ and therefore existence of the inverse $(\nabla^2_{\V z\V z}W)^{-1}$ used in \eqref{EQ:1st_order}-\eqref{EQ:EM_explicit}. However, there might exist points $(\V z^c, \V u^c)$ where the Hessian is non-invertible. 
\HL{The locus of points on the Lagrangian submanifold for which $\det \nabla^2_{\V z\V z}W|_{(\V z^c, \V u^c)}=0$ is called \textit{Maslov-cycle} and its projection on the $\mathcal U$ space determines what in geometric optics is known as \textit{caustics} \cite{Cardin 2015}.}
In these situation, if the more general max-rank condition \eqref{EQ:general_max_rank} holds, then a different set of variables can be used to parameterize the Lagrangian submanifold, without loss of smoothness \cite{Benenti 2011, Cardin 2015}.

\section{Metric considerations from second order analysis}
\label{sec:metric}
In this section, a second order analysis will be carried out to derive a Hessian matrix with tensorial properties. 
While always symmetric, this tensor is not necessarily positive-definite. Based on biological insights we will propose the use of the \textit{squared-Hessian} as a possible metric for quasi-static processes.

\subsection{Equilibrium Forces}
\label{sec:Equilibrium_Forces}

Given the original distinction between internal states and control variables $(\V z, \V u)$, inherent to our systems of interest, we can also distinguish between \textit{internal forces} $\V f_{int}$ and \textit{control forces} $\V f_{ctrl}$ formally as
\begin{equation}\label{EQ:(fint,fu)}
	(\V f_{int}, \V f_{ctrl}) := (-\nabla_{\V z}W, -\nabla_{\V u}W)  .
%			\in T^*(\mathcal Z\times \mathcal U)
\end{equation}
	
When restricted to the equilibrium manifold EM {\eqref{EQ:EM}},  due to condition \eqref{EQ:equil}, equilibrium internal forces are \textit{identically null}
\begin{equation}\label{EQ:f_int=0}\nonumber
		\V 0  = \nabla_{\V z}W|_\text{EM},
\end{equation}
while equilibrium control forces depend on the local branch (indexed by multiplicity index $m$)
\begin{equation}\label{EQ:f_eq}
	{
		\V f_{ctrl}(\V z_m(\V u), \V u) = -\nabla_{\V u}W|_\text{EM} .
	}
\end{equation}

Geometrically, this corresponds to a Lagrangian submanifold on the {control force-space $T^*\mathbb R^K$:}
\begin{equation}\label{EQ:Lambda_red}\nonumber
	\begin{array}{rrl}
	\Lambda=\{(\V  u, \V  f_{ctrl})~:~ 
	& \V f_{ctrl}&=-\nabla_{\V u}W (\V z^*, \V u), \\
	& 0&=\nabla_{\V z}W(\V z^*, \V u), {\text{ for some } \V z^*\in \mathbb R^N } \}.%\subset  T^*\mathbb R^K
	\end{array}
\end{equation}

In previous section,	we showed how a variation of the control inputs $\delta \V u$ around a given input $\V u^*$ induces a variation of internal state $\delta \V z$ prescribed by \eqref{EQ:1st_order}. From definition \eqref{EQ:(fint,fu)}, one can similarly derive variations of forces	
\begin{equation}\label{EQ:delta_f}
	\left\{
	\begin{array}{rllcllll}		
	- \delta \V f_{int} 
		&= \nabla^2_{\V z\V z}W\delta \V z+ \nabla^2_{\V u\V z}W\delta \V u
		&= \V 0 \\
	-\delta \V f_{ctrl} 
		&= \HL{\nabla^2_{\V z\V u}W} \delta \V z+ \nabla^2_{\V u\V u}W\delta \V u 
		&= \V G_m(\V u) \delta \V u ,
	\end{array}
	\right.
\end{equation}
\HL{noting that the first of \eqref{EQ:delta_f} is identically zero due to the equilibrium condition \eqref{EQ:equi_lambda}. The second can be derived by substituting \eqref{EQ:1st_order} in $\delta \V z$.}

\HL{In eq.\eqref{EQ:delta_f}, we introduced the \textit{control Hessian}, defined as } 
\begin{equation}\label{EQ:Gm(u)}
	\V G_m(\V u) := \V G(\V z_m(\V u), \V u) 
\end{equation}
\HL{where}
\begin{equation}\label{EQ:G(z,u)}
{
	\V G(\V z, \V u) := 
	\nabla_{\V u\V u}^2W - \nabla^2_{\V  z\V u}W(\nabla_{\V z\V z}^2W)^{-1}\nabla^2_{\V u\V z}W .
}
\end{equation}

\REMARK{
One can immediately verify that $\V G^T(\V z, \V u)=\V G(\V z, \V u)$ in \eqref{EQ:G(z,u)} is symmetric as well as   $\V G_m(\V u)$ in \eqref{EQ:Gm(u)}.
}
\HL{
\REMARK{\label{rem:Schur}
$\V G(\V z, \V u)$ in \eqref{EQ:G(z,u)} corresponds to  the Schur's complement of the full Hessian of the potential $W(\V z, \V u)$, in particular of the $2\times 2$ block matrix with $\nabla_{\V u \V u}^2W$ and $\nabla_{\V z \V z}^2W$ as diagonal blocks and the mixed derivatives $\nabla_{\V z \V u}^2W$ and $\nabla_{\V u \V z}^2W$ on the remaining blocks, as also clear from \eqref{EQ:delta_f}.
}
}
\REMARK{
A similar formula appears in Gilmore \cite[(5.3)]{Gilmore 1981} in relation to the Taylor's expansion of a potential function defined in terms of internal variables and controls, similarly to $W(\V z, \V u)$. 
}

By making use of \eqref{EQ:1st_order}, one can in fact show that 
\begin{equation}\label{EQ:delta_Wm*}\nonumber
	\begin{array}{rl}
	{\delta^{(2)}} W ^*_m (\V u) := & {W(\V z_m(\V u) +\delta \V z, u+ \delta \V u)-W(\V z_m(\V u), \V u)}\\
	=& 	-\delta \V u^T \V f_{ctrl}(\V z_m(\V u), \V u) 
	+ \frac{1}{2}\delta \V u^T \V G_m(\V u)\delta \V u 
	+\HL{\mathcal O(|\delta\V u|^3)},
	\end{array}
\end{equation}
where 

	\begin{equation}\label{EQ:Wm*}\nonumber
		W_m^*(\V u):=W(\V z_m(\V u), \V u)
	\end{equation}
is the \textit{reduced potential} \HL{for the $m$-th branch} and represents the total energy at equilibrium.

\REMARK{ 
It can be readily verified that the \textit{control Hessian } $\V G_m(\V u)$ is the Hessian of the reduced potential
	\begin{equation}\nonumber
		\V G_m(\V u)=\nabla_{\V u \V u}W_m^*   .
	\end{equation}
}

\subsection{Squared-Hessian metric from Separation Principle}\label{SEC:Separation_Principle}
The control Hessian $\V G_m(\V u)$ in \eqref{EQ:Gm(u)}  is symmetric by definition
	%\footnote{{One can immediately verify that} $\V G^T(\V z, \V u)=\V G(\V z, \V u)$ in \eqref{EQ:G(z,u)} and therefore also $\V G_m(\V u)$ in \eqref{EQ:Gm(u)}  is symmetric.} 
but not necessarily
positive-definite. This will also be shown later on with a simple analytical example. 
{In what follows, we will argue that the squared-Hessian $\V G_m^2 (\V u)$ is in fact a good candidate to determine a metric.
}

\HL{
From an optimal control perspective, we are considering \textit{processes} $s \mapsto \V u(s)$, parameterized by $s\in[0 \ 1]$, starting at $\V u(0)=\V u_0$ and ending at $\V u(1)=\V u_1$. 
We shall assume that two stable equilibria $(\V z_m(\V u_0), \V u_0)$ and $(\V z_m(\V u_1), \V u_1)$ exist \textit{on the same connected}  portion of a given $m$-th branch. 
Under \textit{quasi-static assumptions}, we shall also assume that  such  processes stays on the equilibrium manifold, i.e. $(\V z(u),\V u(s)) \in \text{EM}$ for all $s\in[0\ 1]$.
We would like to single out an input $s\mapsto \V u^*(s)$ which is optimal in some sense, yet to be defined.
}

\HL{To drive the system with input $\V u(s)$,}  the agent will need to  exert a force \HL{$\V f_{ctrl}(\V z_m(\V u(s), \V u(s)))$, in accordance with \eqref{EQ:f_eq}}. 
%and that system will be able to settle at some equilibrium configuration $(\V z_m(\V u^*), \V u^*)$, here assumed stable, 
%{with equilibrium force given by \eqref{EQ:f_eq}.}
%A necessary condition for this to happen is that $\V f_{ctrl}=\V f^{eq}_m(\V u^*)$. 
\HL{Away from criticality, a `small' change in the input $\V u^*+\delta \V u$ will result in a small variation in control force $\delta \V f_{ctrl}$.  
}
Next, we shall study the relationship between these two variations.

A first idea could be trying to minimize the 
	\textit{total effort} $\|\V f_{ctrl}(s)\|^2$.
However, \HL{at any `time' $s$}, 
the amount of force $ \V f_{ctrl}(\V z_m(\V u\HL{(s)}),\V u\HL{(s)})$ given by \eqref{EQ:f_eq} is necessary to maintain equilibrium and, in many applications, it is likely to represent the bulk of the total effort.   
A biological insight, known as \textit{separation principle}, suggests that the equilibrium force {$\V f_{ctrl}(\V z_m(\V u\HL{(s)}),\V u\HL{(s)})$} should be in fact \textit{factored out} of the optimization process as it represents a necessary effort to maintain equilibrium. 
In other words, a biological agent is more likely to minimize the remaining \textit{residual effort} 
\HL{
$$\delta \V f_{ctrl}(\V u(s)) := \V f_{ctrl}(\V z_m(\V u(s)),\V u(s))- \V f_{ctrl}(\V z_m(\V u^*(s)),\V u^*(s))$$ 
}
along the process $\V u(s)$. At infinitesimal level, the residual effort is represented by $\delta \V f_{ctrl}$ in \eqref{EQ:delta_f} and its minimization corresponds to minimizing

\begin{equation}\nonumber
	\begin{matrix}
	\|\delta \V f_{ctrl}\|^2
		&= \delta \V f_{ctrl}^T \delta \V f_{ctrl} \hfill\\
		&= \delta \V u^T \V G_m^T(\V u)\V G_m(\V u)\delta \V u \hfill\\
		&= { \delta \V u^T \V G_m^2(\V u)\delta \V u .} \hfill
	\end{matrix}
\end{equation}

\REMARK{  By symmetry of the Hessians, $\V G_m^T(\V u) \V G_m(\V u) = \V G_m^2(\V u)$, the squared-Hessian is \textit{symmetric and non-negative}, by definition. 
}
As shown in \HL{Appendix \ref{sec:Appendix_Hessian}, in particular in Eq.~\eqref{eq:Hessian_covariance} and Eq.~\eqref{eq:squared_Hessian_covariance}}, both $\V G_m$ and $\V G_m^2$ behave tensorially and the metric $\V G_m^2$ is simply a sort of (vertical) \textit{pull-back}, \HL{see Eq.~\eqref{eq:Vertical_pullback} in Appendix \ref{sec:Appendix_squaredHessian},} on the Lagrangian submanifold $\Lambda_m$ of any given natural metric on \HL{the} control manifold $\mathbb R^K$.
We therefore propose to use the squared-Hessian $\V G_m^2(\V u)$ as a metric to define optimal processes as paths $\V u(s):[0\ 1]\rightarrow\mathbb R^K$ extremizing the following functional 
\begin{equation}\label{EQ:J_functional}
	J:=\int_0^1 {\sqrt{\HL{\Vdot u^T} \V G^2_m \HL{\Vdot u} } }\,ds   ,
\end{equation}
where   $\HL{\Vdot u}:=d\V u/ds$.

\REMARK{
\HL{The structure of the integrand in \eqref{EQ:J_functional}, containing $\Vdot u$, implies that we seek to minimize variations of $\V f_{ctrl}$ both along $\V u^*(s)$ and along  curves $\V u(s)$ homotopic to $\V u^*(s)$.}
The square root is essential in order to put in evidence that the
solution has to be independent of the parametrization. 
}
This is not a great analytical bias: the relations between the stationary solutions of $J = \int_0^1\sqrt{\HL{\Vdot u}^T \V G^2_m \HL{\Vdot u} }ds$ and $\bar J = \int_0^1{\HL{\Vdot u}^T \V G^2_m\HL{\Vdot u   } }ds$
are well known:  
{ 
if $\V u(s)$ stationarizes $\bar J$  then it stationarizes $J$ as well. In the other case, if $\V u(s)$ stationarizes $J$ then it can be shown that there exists a suitable reparameterization such that the new curve stationarizes $\bar J$ as well.
}

\subsection{Gauss' Principle and the $\V G^2$ metric}
%Gauss' Principle \cite{Gauss 1829} is an ancestral tool for describing dynamics of  ideal, constrained mechanical system. Born in a genuine non-static environment, it states that along the admissible motions, point by point in the host tangent space, \textit{the square of the norm of the global reaction forces has to be minimum for any admissible variations of the accelerations}. This notable Principle, described in an intrinsic geometrical framework \cite{Cardin and Zanotto 1989}, has been often utilized in many contexts: statistical mechanics, see \cite[p. 320]{Gavallotti 1999}, and fluid dynamics,  \cite[p. 444, p. 451]{Gavallotti 2013}; furthermore, by thinking of our aims, in various robotic settings, e.g. \cite{Lilov and Lorer 1982}.

Here, we try to recognize inside our developed point of view the natural
germ of that powerful Gaussian idea, even though in a quasi-static context.
The starting key point is to observe that we arrived to declare the minimization
of a strict analogue of the reaction forces, generated by the acting
control parameters $\V u(s)$, as described in \HL{a} previous section.

In our control context, we encounter a more rich scenario. First \HL{of} all,
we must build the Lagrangian submanifold of the equilibria $\Lambda \subset T^*\mathbb R^K$. Then, we compute also the parametrizations of the transversal branches
of $\Lambda$:
$$z_1(\V u), \dots, z_m(\V u), \dots , z_M(\V u).$$

There is a radical difference between standard constrained mechanics and
the actual controlled dynamics: if we are at $\bar{\V u}$, over the {$m$-th} branch, a
reasonable imported \HL{Gauss'} Principle says us that the right direction $\delta \V u$ to move from $\bar {\V u}$ to $\bar {\V u}+\delta \V u$ is in correspondence to the `minimum' (normalized) eigenvector $\delta \V u$ of $\V G^2$. 
However, our control problem is a bit different: we want to move between extreme points, from  $\V u_0$ to $\V u_1$, relative to two pre-defined equilibrium configurations {$\V z_{m0}(\V u_0)$ and $\V z_{m1}(\V u_1)$}

 Necessarily, we are brought back to a non-local variational principle (unlike \HL{Gauss'}). In this order of ideas, this corresponds to taking  the minimum of the integral of $|| \V f_{ctrl}||^2$, denoting the square of the norm of the ``reaction forces''.
\HL{The control problem will therefore be to determine}
 a multi-branch trajectory %
\begin{equation}\nonumber
	[0\ 1]\ni s \mapsto (\V u(s), m(s)) \in \mathbb R^K \times \{1,  \dots, M\}
\end{equation}
among the curves with $\V u(0)=\V u_0$ and $\V u(1)=\V u_1$ and such that
\begin{equation}\nonumber
	\inf \int_0^1 \HL{\Vdot u(s)^T \V G^2_{m(s)}\HL{\Vdot u(s)} } \,ds   .
\end{equation}
\REMARK{
Solutions to this variational principle are allowed to \textit{change branch} along their quasi-static evolution, e.g. at some $s^*$, jumping from $\V z_{m1}(s^*)$ to $\V z_{m2}(s^*)$, where
	$$\lim_{s \to (s^*)^-} m(s) = {m_1}, 
		\hspace{2cm} 
	\lim_{s \to (s^*)^+} m(s) = {m_2}.$$
}

\section{Toy Example: manipulation of an elastically driven inverted pendulum}
\label{sec:toy_model}
In this section, the theory developed \HL{in the first part of the paper} will be applied to a simple 2D pendulum elastically driven by an agent. 
\HL{
    The basic mechanical definition of the problem is given in Sec.~\ref{sec:inverted_pendulum_definition}. In this short subsection, all one needs to define is $i)$ the \textit{kinematics} of the problem, including the transformations between the various frames of reference, as commonly done an any robotics problem; $ii)$ the \textit{energetics} of the problem, i.e. the scalar potential $W(\V z, \V u)$ which typically consists of the sum of all gravitational potentials (one term per rigid body) and elastic potentials (one term for each generalized spring between bodies).  The subtle step involved in the latter point is the choice of variables, i.e. deciding which ones are to represent the internal states and which ones are to be controlled by the agent. We believe the Engineering of the problem  lies in this very choice.

After the problem definition, in Sec.~ \ref{sec:ToyModel_constant_spring} we consider the case of an agent interacting with an inverted pendulum via a linear and constant spring. In this very simple case, all calculations could   easily be derived by hand and the optimal control problem can be solved without any approximation. The main purpose of the case is to illustrate how simply the theory described in the first part of the paper can be applied to solve optimal manipulation problems, at least while remaining on a single branch of the equilibrium manifold.
   
To address more realistic scenarios, in Sec.~\ref{sec:ToyModel_numerical} we show how mechanical interaction can be approximated via a smooth albeit high-nonlinear spring. This is akin to regularizing the mechanical contact which is intrinsically a non-smooth problem \cite{Brogliato}. The smoothness of the problem allows for standard numerical approaches to be applied. We conclude by showing how Graph Theory can be used to provide a connected view of the set of discrete numerical solutions. The connectivity of the constructed graph includes weights derived from the proposed squared-Hessian metric. This allows the manipulation problem to be rephrased and numerically  solved as minimum path of a weighted graph.

}
\subsection{Elastically driven inverted pendulum}
\label{sec:inverted_pendulum_definition}
\begin{figure}[h]
	\centering
	\includegraphics[scale=0.9]{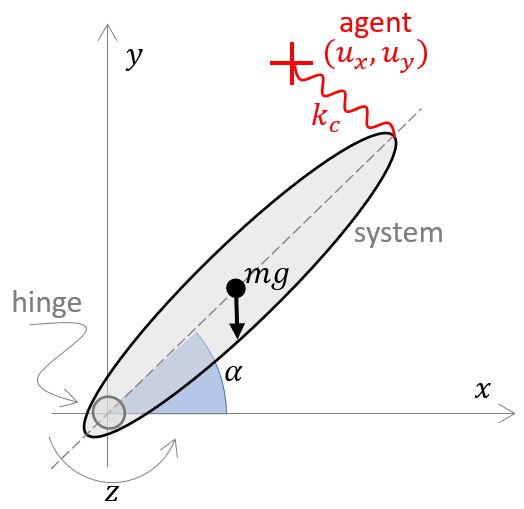}
	\caption{Inverted Pendulum}
	\label{fig:Ex1}
\end{figure}

With reference to Fig.~\ref{fig:Ex1}, we consider the task of stabilizing an inverted pendulum in its upward position, in the space frame $\{0\}$. 
The pendulum consists of a rigid body $\mathcal B\subset\mathbb E^2$ (i.e. {the plane} with the usual Euclidean distance). For simplicity, we shall assume a uniform, rectangular bar of length $L_0$ with center of mass ($com$) located at its geometric center.

We shall consider a moving-frame $\{1\}$ attached to the $com$ of the pendulum and rotating with it by angle $\alpha$ (see Fig.\ref{fig:Ex1}) with respect to the space-frame $\{0\}$. For convenience, we shall define the (unit-length) axis $\V n_\alpha$ directed along the major axis,  the (unit-length) axis $\V n_\alpha^\perp$ perpendicular to it as well as the rotation matrix $\V R_\alpha$  as follows:

\begin{equation}\label{eq:SE2_transformations}
	\V n_\alpha = \MAT{\cos\alpha \\ \sin\alpha} \hspace{.6cm}
	\V n_\alpha^\perp = \MAT{-\sin\alpha \\ \cos\alpha} \hspace{.6cm}
	\V R_\alpha = \MAT{\cos\alpha & -\sin\alpha\\ \sin\alpha & \cos\alpha} .
\end{equation}

As shown in Fig.\ref{fig:Ex1}, the pendulum is subject \HL{to} gravitational potential $W_{grav}$ (due to a constant force $mg$, pointing downwards and acting at $com$) as well as to the elastic (control) potential $W_{ctrl}$ due to a linear spring of  stiffness $k_c$ attached between an agent-controlled position $\V u=(u_x, u_y)$ and the tip of the pendulum $\V c = L_0\V n_\alpha$.
The total potential $W$ can then be expressed as 
\begin{equation}\label{EQ:Ex1:W}
	\begin{matrix}
		W(\alpha, \V u) &=	&  W_{grav} + W_{ctrl} \hfill\\
		&= & \frac{1}{2}mgL_0\sin\alpha + \frac{1}{2}k_c\|\V u-L_0\V n_\alpha\|^2 \hfill\\
		&= & \frac{1}{2}mgL_0\sin\alpha + \frac{1}{2}k_c((u_x-L_0\cos\alpha)^2 + (u_y-L_0\sin\alpha)^2) . \hfill
	\end{matrix}
\end{equation}

\HL{
Note that the definition of the kinematics , i.e. \eqref{eq:SE2_transformations}, and the energetics of the problem, i.e. \eqref{EQ:Ex1:W}, are the main two steps requiring the user to `manually' input data from the problem. The remaining steps involve  algebraic manipulation which can, in principle, be automated with any Computer Algebra System. 
As an example, all symbolic computations and numerical results presented in this section are derived via MATLAB \textit{livescripts} which are provided as supplementary material, along with their printouts. 
}

\HL{ 
\subsection{Linear spring case: an analytical study}  
\label{sec:ToyModel_constant_spring}
}

As a first case, we shall consider a linear spring of constant stiffness $k_c$. In this case, first and higher derivatives of the potential energy \eqref{EQ:Ex1:W} become straightforward and can be tackled analytically.

\subsubsection{Equilibria}
For a given control input $\V u$, the equilibrium position(s) of the pendulum can \HL{be} found from \eqref{EQ:equil}, where the variable $\V z$ is identified with the scalar $\alpha$ (the pendulum angle) in this example
\begin{equation}\label{EQ:Ex1_Wz=0}
	\begin{matrix}
		0&=& \nabla_{\V z}W \hfill\\
		 &=& k_cL_0(u_x\sin\alpha  -(u_y-mg/2k_c)\cos\alpha) \\
		 &=& -k_cL_0  (\V n_\alpha^\perp)^T(\V u-\V u_{crit}) ,\hfill
	\end{matrix}
\end{equation}
where $\V u_{crit}:=[0\ mg/2k_c]^T$ denotes a \textit{critical point} where the system will display instability, as it will be shown soon.

As the pendulum is a physical system in the standard Euclidean space (2D, for simplicity), what follows will be make use of the familiar Euclidean norm.
It is apparent that \eqref{EQ:Ex1_Wz=0} expresses an orthogonality condition between the vector $\V u-\V u_{crit}$ and the unit vector $\V n_\alpha^\perp$ (which is, in turn, perpendicular to the pendulum). 
Therefore $\V u - \V u_{crit}$ is either aligned (+) or anti-aligned (-) with the pendulum direction $\V n_\alpha$, i.e. \eqref{EQ:Ex1_Wz=0} is equivalent to

\begin{equation}\label{EQ:Ex1_LuNa}
	\V u-\V u_{crit}	= \pm L_u \V n_\alpha ,
\end{equation}
where, as also shown in Fig.\ref{fig:Ex1}, \textit{control length} $L_u$ is defined as 
%$L_u$ represents the magnitude $\|\V u-\V u_{crit}\|$, i.e.
\begin{equation}\label{EQ:L_u}
	L_u:= \|\V u-\V u_{crit}\| = \sqrt{u_x^2+(u_y-mg/2k_c)^2} .
\end{equation}

Away from the critical point (i.e. $\V u \neq \V u_{crit}$ or, equivalently, $L_u>0$), analytical solutions of \eqref{EQ:Ex1_Wz=0} can be described as 
\begin{equation}\label{EQ:Ex1_atan2}
	\V z_m^*(\V u) := \text{atan2}(u_y-mg/2k_c, u_x) +m\pi,
\end{equation}
where $m=0$ corresponds to the aligned (+) case and $m=1$ corresponds to the anti-aligned case (-). Given the periodicity of the problem, any other value of multeplicity $m$ can be considered equivalent to $m=0$ (if even) or to $m=1$ (if odd). Below, it will be shown that $m=0$ is a stable solution while $m=1$ is unstable.
We shall therefore only consider the stable solution
\begin{equation}\label{EQ:Ex1:alpha*}
	\alpha^*:=\V z_0^*(\V u)=\text{atan2}(u_y-mg/2k_c, u_x).
\end{equation}
Restricting the potential to the equilibrium manifold leads to the \textit{reduced potential} 
\begin{equation}\label{EQ:Ex1:W*}
	W^*(\V u):=W(\alpha^*, \V u) = \frac{1}{2}k_c(L_u - L_0)^2 + \frac{1}{2}mg\, u_y +const.
\end{equation}
This can be verified by noting that the stable equilibrium condition \eqref{EQ:Ex1_LuNa}, i.e. the one with the positive sign, corresponds to $u_x=L_u\cos\alpha$ and $u_y=L_u\sin\alpha+mg/2l_c$ which, substituted back in \eqref{EQ:Ex1:W}, after some algebraic manipulation, lead to \HL{Eq.}\eqref{EQ:Ex1:W*}.

\REMARK{
	The constant and linear terms in \eqref{EQ:Ex1:W*} are inessential for the purpose of evaluating the Hessian which, therefore, only depends on the quadratic term $1/2k_c^2(L_u-L_0)^2$. Such term captures the \textit{radial push/pull}, it only affects the reaction forces at the hinge, not the ones against gravity.
}

\subsubsection{Stability}
To analyze stability, we will need to evaluate higher derivatives with respect to the parameter $\alpha$. To this end, it is convenient to note the following identities
\begin{equation}\label{EQ:Ex1_identities}\nonumber
	\frac{\partial}{\partial \alpha} \V n_\alpha = \V n_\alpha^\perp,
	\hspace{1cm}
	\frac{\partial}{\partial \alpha} \V n_\alpha^\perp = -\V n_\alpha,
\end{equation}
through which one immediately evaluates the Hessian

$$
	\begin{matrix}
		\nabla_{\V z\V z}^2W &= & -k_cL_0  (\frac{\partial}{\partial \alpha}\V n_\alpha^\perp)^T(\V u-\V u_{crit}) \hfill \\
		&=&+k_cL_0\V n_\alpha^T (\V u-\V u_{crit}) \hfill \\
		&=&\pm k_cL_0\V n_\alpha^T L_u\V n_\alpha \hfill \\
		&=&\pm k_cL_0 L_u  . \hfill \\
	\end{matrix}
$$

While $k_c>0$, $L_0>0$ are strictly positive constants, the term $L_u$ defined in \HL{Eq.}\eqref{EQ:L_u} is strictly positive if $\V u\neq \V u_{crit}$. Therefore, only the aligned (+) solution in \eqref{EQ:Ex1_LuNa} ensures stability via a strictly positive Hessian, therefore we shall only consider the non-negative case

\begin{equation}\label{EQ:Ex1_Wzz}\nonumber
	\nabla_{\V z\V z}^2W = k_cL_0 L_u ,
\end{equation}
which is in fact positive definite as long as $\V u\neq \V u_{crit}$.
\subsubsection{Second Derivatives}
%In addition to the Hessian \eqref{EQ:Ex1_Wzz}, 
In order to evaluate the reduced Hessian $\V G_m(\V u)$ as in \eqref{EQ:Gm(u)}, the remaining second derivatives have to be evaluated. Differentiating the total potential $W$ in \eqref{EQ:Ex1:W} twice with respect to $\V u$ immediately leads to
\begin{equation}\nonumber
	\nabla_{\V u\V u}W=\MAT{k_c & 0\\0 & k_c}   .
\end{equation}
Differentiating $\nabla_{\V z}W$ in \eqref{EQ:Ex1_Wz=0} with respect to $\V u$  returns 
\begin{equation}\nonumber
	\nabla_{\V z\V u}W=k_cL_0 \MAT{\sin\alpha\\ -\cos\alpha}= - k_cL_0\V n_\alpha^\perp
\end{equation}
and, in general,  $\nabla_{\V u\V z}W=\nabla_{\V z\V u}^TW$.

Therefore the matrix $\V G(\V z, \V u)$ in \eqref{EQ:G(z,u)} can be rewritten as
\begin{equation}\label{EQ:Ex1_G}\nonumber
	\begin{matrix}
		\V G(\alpha, \V u) &=& \nabla_{\V u\V u}W - \nabla_{\V z\V u}W(\nabla_{\V z\V z}W)^{-1}\nabla_{\V u\V z}W \\
		&=&{\tiny \MAT{k_c&0\\0&k_c}} - k_c\frac{L_0}{L_u}\V n_\alpha^\perp(\V n_\alpha^\perp)^T\\
		&=& \V R_\alpha {\tiny \MAT{k_c&0\\0&k_c\left(1-\frac{L_0}{L_u}\right)}} 
		\V R_\alpha^T,
	\end{matrix}
\end{equation}
while the reduced Hessian \eqref{EQ:Gm(u)} can be derived by evaluating the general Hessian at the equilibrium $\alpha^*$ in \eqref{EQ:Ex1:alpha*}
\begin{equation}\label{EQ:Ex1_G*}\nonumber
	\V 	G_m(\V u):= \V G(\alpha^*, \V u)
	 =  \V R_{\alpha^*} {\tiny \MAT{k_c&0\\0&k_c\left(1-\frac{L_0}{L_u}\right)} }
	 \V R_{\alpha^*}^T.
\end{equation}

From this form, one can immediately verify that $\V n_{\alpha^*}$ and $\V n_{\alpha^*}^\perp$ are the two \textit{eigenvectors} corresponding, respectively, to a constant and strictly positive eigenvalue $k_c$ and to a second eigenvalue $k_c(1-L_0/L_u)$ which is negative when the pendulum is compressed ($L_u<L_0$), zero when $L_u=L_0$ and positive otherwise.

On the other hand, by definition, the squared-Hessian $\V G_m^2(\V u)$ ensures non-negative eigenvalues
\begin{equation}\label{EQ:Ex1_G*^2}\nonumber
	\V 	G_m^2(\V u)
	=  \V R_{\alpha^*} {\tiny \MAT{k_c^2&0\\0&k_c^2\left(1-\frac{L_0}{L_u}\right)^2}} \V R_{\alpha^*}^T.
\end{equation}
\REMARK{
The \textit{eigenvectors are unaltered} by this operation while the eigenvalues are squared.}

\subsubsection{Optimal Control via Squared-Hessian}

For the (simple) problem considered in this section, the equilibrium manifold is a \textit{ruled surface} defined by \eqref{EQ:Ex1_atan2} in the 3-parameter space $(\alpha, u_x, u_y)$: for every value of $\alpha$, the locus of solutions $(u_x, u_y)$ is a line going through $\V u$ and $\V u_{crit}$ (we shall always assume $\V u\neq \V u_{crit}$).

As previously mentioned, the equilibrium conditions \eqref{EQ:Ex1_atan2} and \eqref{EQ:Ex1_LuNa} are equivalent but the latter allows treating $\alpha$ as an independent parameter, {replacing therefore `$s$', for this specific example}. This is useful when considering problems such as driving the pendulum quasi-statically (i.e. through consecutive equilibria) from an initial position $\alpha_1=-\pi/2$ to a final up-right position $\alpha_2=+\pi/2$ of the pendulum. Specifically, the stable condition (+) \eqref{EQ:Ex1_LuNa} can be rewritten as 
\begin{equation}\nonumber
	\left\{\begin{matrix}
		u_x(\alpha) = L_u(\alpha) \cos\alpha \hfill \\
		u_y(\alpha) = L_u (\alpha) \sin\alpha+mg/2k_c ,
	\end{matrix}
	\right.
\end{equation}
where the control length $L_u$ is now considered as a dependent variable.
The optimal problem can now be formulated as a minimization of the action
\begin{equation}\label{EQ:Ex1_J}
	J = \int_{\alpha_1}^{\alpha_2}{\V u'}^T \V G_m^2 {\V u'}\, d\alpha ,
\end{equation}
with initial and final conditions
\begin{equation}\label{EQ:Ex1:BoundaryContitions}\nonumber
		\left\{\begin{matrix}
		L_u(\alpha_1) = L_1 \\
		L_u(\alpha_2) = L_2 . \\
	\end{matrix}
	\right.
\end{equation}

\HL{Note: in this section, since the time variable $s$ is now replaced with the scalar angle $\alpha$, we shall use} the `prime' $(')$ operator to denote differentiation with respect to $\alpha$ and, from \eqref{EQ:Ex1_LuNa}, one derives
\begin{equation}\nonumber
	{\V u'}=(L_u \V n_\alpha)'
		=L_u' \V n_\alpha + L_u \V n_\alpha^\perp  .
\end{equation}
Using the identities $\V R_\alpha^T \V n_\alpha = [1\ 0]^T$ and  $\V R_\alpha^T \V n_\alpha^\perp = [0\ 1]^T$, the integrand in \eqref{EQ:Ex1_J} corresponds to the (`pendulum') Lagrangian \HL{$\mathcal L^{pend}$} of the optimization problem and can be expanded as

\begin{equation}\label{EQ:Ex1_Lagrangian}\nonumber
	\begin{matrix}
		{\V u'}^T \V G_m^2 {\V u'} &=& 
			(L_u' \V n_\alpha^T + L_u (\V n_\alpha^\perp)^T)
			\V R_{\alpha} {\tiny\MAT{k_c^2&0\\0&k_c^2\left(1-\frac{L_0}{L_u}\right)^2} }
			\V R_{\alpha}^T
			(L_u' \V n_\alpha + L_u \V n_\alpha^\perp) \hfill\\
		&=& 
			(L_u' [1\ 0] + L_u [0\ 1])
		 {\tiny\MAT{k_c^2&0\\0&k_c^2\left(1-\frac{L_0}{L_u}\right)^2} }
			(L_u' [1\ 0]^T + L_u [0\ 1]^T) \hfill \\ 
		&=& k_c^2 
			{L_u'}^2 +k_c^2 L_u^2(1-\frac{L_0}{L_u})^2 \hfill \\
		&=& k_c^2 
			\underbrace{({L_u'}^2 +L_u^2 + L_0^2 - 2L_0L_u)}_{\HL{\mathcal L^{pend}}} . \hfill \\
	\end{matrix}
\end{equation}

An optimal solution to minimization of the action integral \eqref{EQ:Ex1_J} with Lagrangian \HL{$\mathcal L^{pend}$} can be found from Euler Lagrange equations

\begin{equation}\nonumber
	\left(\frac{\partial \HL{\mathcal L^{pend}}}{\partial L_u'}\right)'-\frac{\partial \HL{\mathcal L^{pend}}}{\partial L_u}=0
	\hspace{5mm}\equiv \hspace{5mm}
	L_u'' = L_u-L_0  .
\end{equation}

\REMARK{
This is a second order, linear equation which does not depend on parameters such as spring stiffness $k_c$ or pendulum weight $mg$. 
}
To make it independent of the pendulum length as well, we can simply consider the \textit{normalized length} 
\begin{equation}\nonumber
	\lambda_u(\alpha):=L_u(\alpha)/L_0
\end{equation}
and therefore the normalized differential equation
\begin{equation}\nonumber
	\lambda_u'' = \lambda_u - 1,
\end{equation}
with  solutions
\begin{equation}\nonumber
	\lambda_u(\alpha) = a e^{+\alpha} + b e^{-\alpha} + 1.
\end{equation}

Considering normalized initial and final lengths $\lambda_1:=L_1/L_0$ and $\lambda_2:=L_2/L_0$, the coefficient $a$ and $b$ can be determined as 
\begin{equation}\nonumber
	\MAT{a\\ b} = \MAT{e^{+\alpha_1} & e^{-\alpha_1} \\ e^{+\alpha_2}&e^{-\alpha_2}}^{-1}
	\MAT{\lambda_1-1 \\ \lambda_2-1},
\end{equation}
which gives the closed-form solution:
\begin{equation}\nonumber
	\lambda_u(\alpha) = \MAT{e^{+\alpha}& e^{-\alpha}} \MAT{e^{+\alpha_1} & e^{-\alpha_1} \\ e^{+\alpha_2}&e^{-\alpha_2}}^{-1}
	\MAT{\lambda_1-1 \\ \lambda_2-1} +1.
\end{equation}
Note that the 2$\times$2 matrix is invertible \textit{iff} $\alpha_1\neq \alpha_2$.

\begin{figure}[h!]
	\centering
	\includegraphics[width=0.7\linewidth]{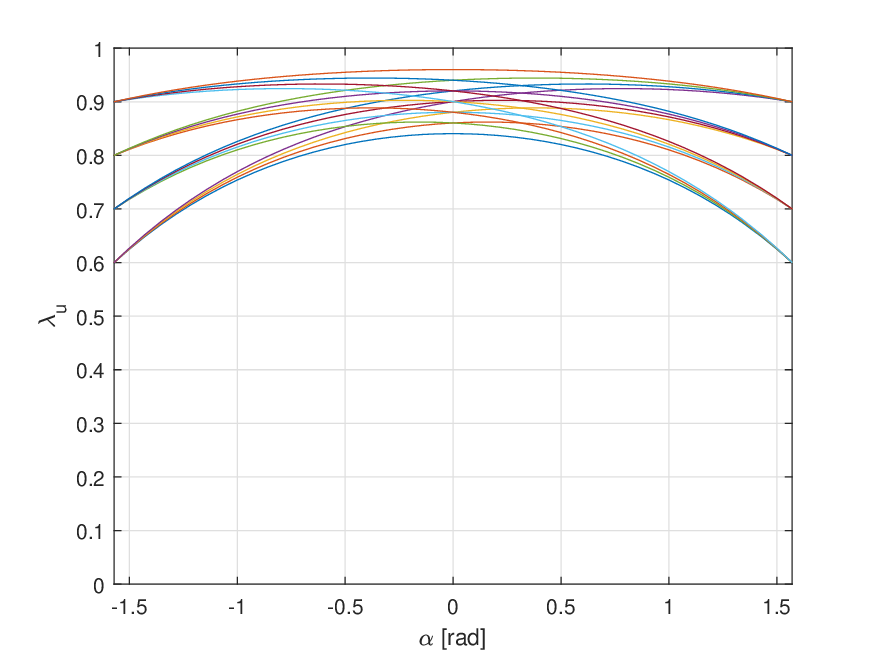} 
	\includegraphics[width=0.7\linewidth]{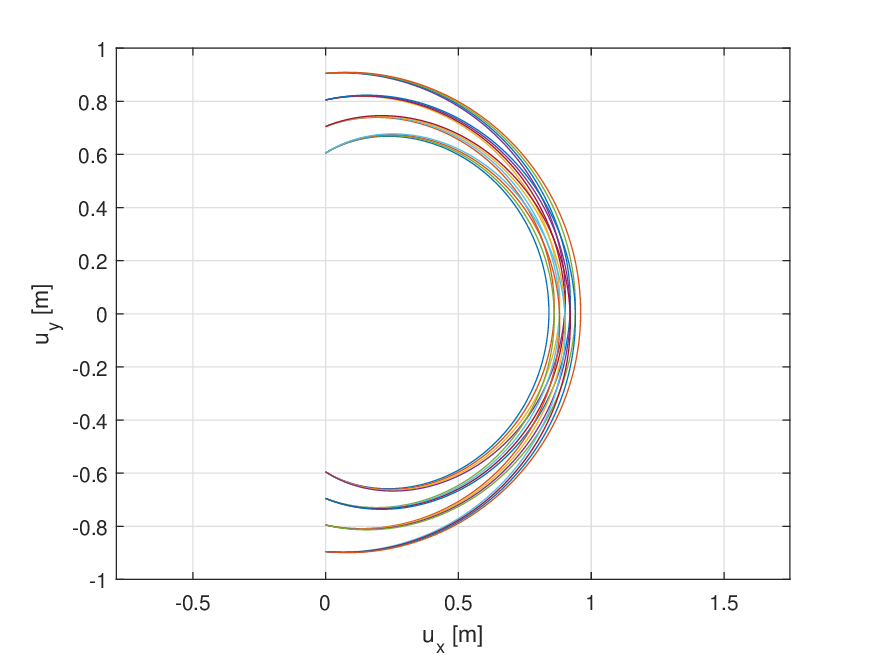} \\
	\caption{Optimal curves $\lambda_u(\alpha)$ vs. $\alpha$ (top) and control positions \HL{$u_y(\alpha)$ vs. $u_x(\alpha)$ } for different values of initial and final length}
\end{figure}

\subsection{Nonlinear spring case: towards robotic manipulation, a numerical study}
\label{sec:ToyModel_numerical}
\begin{figure}[h]
	\centering
	\includegraphics[scale = 0.9]{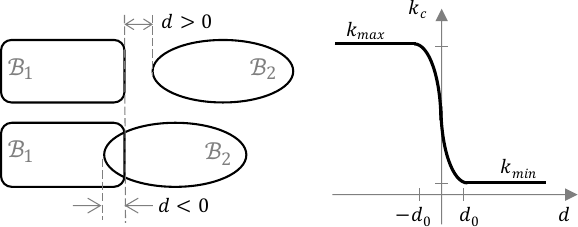} 
	\caption{Regularization of mechanical contact: based on the penetration variable $d$, a spring with nonlinear stiffness $k_c$ is designed assume high-stiffness levels ($k_{max}$) when the bodies are in contact ($d<0$) and low-stiffness levels ($k_{min}$) when not in contact $d>0$. A steep yet smooth transition occurs in an interval $[-d_0\ d_0]$.}
	\label{fig:Ex2_Kc}
\end{figure}

The example in Fig.~\ref{fig:Ex1} was devised to be simple enough to allow for exact, analytical solutions. Most of real-life scenarios, however, do not share this luxury. Specifically for robotic manipulation, one major hurdle is the non-smooth nature of mechanical contact \cite{Brogliato} for which there exist various possible approaches. Here we would like to be able to make use of the optimal control framework derived so far, for which smoothness is essential. We therefore propose to \textit{regularize} contact between two bodies (e.g. $\mathcal B_1$ and $\mathcal B_2$) via a non-linear stiffness $k_c$ which depends on the \textit{inter-penetration} $d$, as sketched  in Fig.~\ref{fig:Ex2_Kc}.

Inter-penetration is a \textit{signed distance}: 
if two objects are in mechanical contact, we allow for some degree of interpenetration ($d<0$) and high contact forces arise from the assumption of a high level of stiffness ($k_{max}$). On the other hand, when the two objects are not in contact ($d>0$), one would expect no forces and this is modeled via a very low level of stiffness ($k_{min}\ll k_{max}$). 
\REMARK{
A very low but non-zero stiffness is \textit{essential} as it provides useful non-zero gradients for the (virtual) agent to find a path even when not in contact. It turns out to be a very useful expedient to make robotic agents aware of their surroundings. For example, in \cite{Kana et al (2021)}, virtual `proxies' are introduced to guide the robot towards the goal before any contact is even made.
} 

We then  assume a steep yet smooth transition (around $d\approx 0$) as sketched in the graph in Fig.~\ref{fig:Ex2_Kc}. Analytically, this can be modeled via smooth functions such as 

\begin{equation}\label{EQ:k_c}
	k_c(d)=k_{min} +   \frac{ 1-\tanh(d/d_0) }{2} k_{max}  ,
\end{equation}
where $d_0$ is a parameter capturing the penetration depth.

Computing the inter-penetration of rigid bodies for generic shapes can be computationally expensive and  is typically done via numerical methods \cite{Kurtz and Lin 2022}. The task is much simpler when one of the bodies can be assumed to be a point, as in the case of the agent in Fig.\ref{fig:Ex1}, represented by a point of coordinates $\V u=(u_x, u_y)$.
\HL{
	In robotics, this is not uncommon. In the so-called `force control mode', robots are programmed to impart actual forces based on virtual springs, such as $k_c$, virtual points such as  $\V u$ and physical points, such as  $\V c$, whose coordinates can be measured, for example, via vision systems.
} 

With reference to the specific 2D problem in Fig~\ref{fig:Ex1}, to compute the penetration $d$ by the agent, localized at given point $\V u$, within the rigid body of the pendulum, a geometrical definition of the shape of the pendulum is required. Both in 2D as well as in 3D, a simple way to proceed analytically is to model the rigid body $\mathcal B$ as a \textit{super-ellipse} \cite{Jaklic et al 2000}, for which the pseudo-distance $d$ between agent and body $\mathcal B$ is simply given as 
	$$d:=\Delta(u_x, u_y), $$
where $\Delta:\mathbb E^2\to \mathbb R$  is the super-ellipse \textit{inside-outside function}
\begin{equation}\label{EQ:Delta}
	\Delta: (x,y)\mapsto \left(\frac{x}{a}\right)^{2/\epsilon} + \left(\frac{y}{b}\right)^{2/\epsilon}-1
\end{equation}
and where $0<\epsilon<2$ defines the shape of the super-ellipse ($\epsilon=1$ for an ellipse, $\epsilon<1$ for a super-ellipse as in this work), while $a$ and $b$ represent half-length and half-width of the body $\mathcal B$, respectively.

\subsubsection{Equilibrium Manifold}
The configuration space $Q$ of the agent+pendulum can be parameterized by the (scalar) rotation $z\equiv \alpha \in \mathbb R$ of the pendulum and the 2D position $\V u\equiv (u_x, u_y)\in \mathbb R^2$ of the agent, resulting in the three-dimensional manifold:
\begin{equation}\nonumber
	Q:= S^1\times \mathbb R^2 .
\end{equation}
As our analysis is local, 
{ we shall assume  working in an open \HL{set} of $\mathbb R\times \mathbb R^2$.	

The total potential $W$ will comprise a gravitational potential $W_{grav}$ and  the (control) elastic potential $W_{ctrl}$ and will be formally similar to \eqref{EQ:Ex1:W}.
The only difference is due to the nonlinearity of the stiffness $k_c$, now function of the inter-penetration $d$ which can readily be evaluated via the super-ellipse function $\Delta(\cdot)$ defined in \eqref{EQ:Delta}. This valuation is more conveniently done in the moving frame: 
\begin{equation}\label{EQ:W_c}\nonumber
	W_{ctrl} = \frac{1}{2} k_c(\Delta(\tilde {\V u})) \ \|\tilde {\V u} - \tilde {\V c}\|^2,
\end{equation}
where $\|\cdot\|^2$ is the Euclidean distance; the tilde $\tilde \cdot$ denotes coordinates expressed in the moving frame $\{1\}$, in which, $\tilde{\V c}$ is a constant vector and $\tilde {\V u}:= \V T_1^{-1}\circ \V u$ .

Thanks to the smoothness of the functions \eqref{EQ:k_c} and \eqref{EQ:Delta}, 
gradients $\nabla_{\V z}W$, $\nabla_{\V u}W$ as well as the second derivatives $\nabla^2_{\V u\V z}W$, $\nabla^2_{\V z\V u}W$, and the Hessian $\nabla^2_{\V z}W$
can be computed analytically in symbolic form, via the MATLAB Symbolic Toolbox or any other Computer Algebra System.

Using the numerical values given in  Table~\ref{TAB:values}, a numerical solver ({\tt fsolve} in MATLAB) was used to solve the equilibrium equation \eqref{EQ:equil} for an array of inputs $u_x$ and $u_y$, each spanning across 150\% of workspace ($\pm L_0$). 

\begin{table}[h!]\centering	
	\begin{tabular}{c|ccl}
		\textit{parameter} & \textit{value} & \textit{units} & \textit{description}\\
		\hline
		$L_0$ & 1 & $m$ & pendulum length \\
		$W_0$ & 0.1 & $m$ & pendulum width \\
		$mg$ & $10$ &N & weight force\\
		$k_{min}$ & $1$ & N/m & stiffness \\
		$k_{max}$ & $10^4$& N/m & stiffness \\
		$\epsilon$ & 0.1 &-& super-ellipse factor
	\end{tabular}
	\caption{Values used in the numerical simulations.}
	\label{TAB:values}
\end{table}

\begin{figure}[H]
	\centering
	\includegraphics[width=0.70\linewidth ]{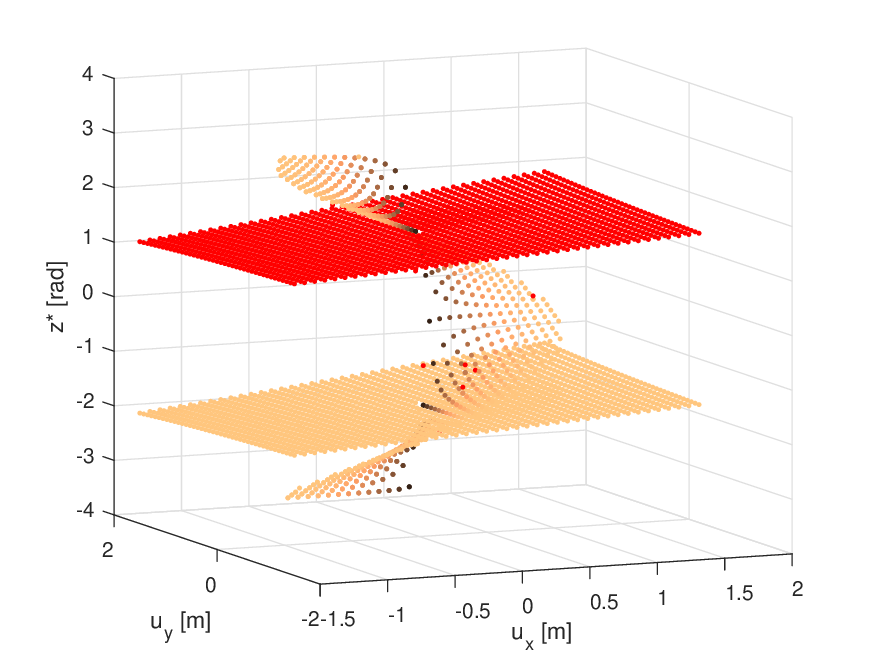}
	\\ 	a) \\
	\includegraphics[width=0.70\linewidth]{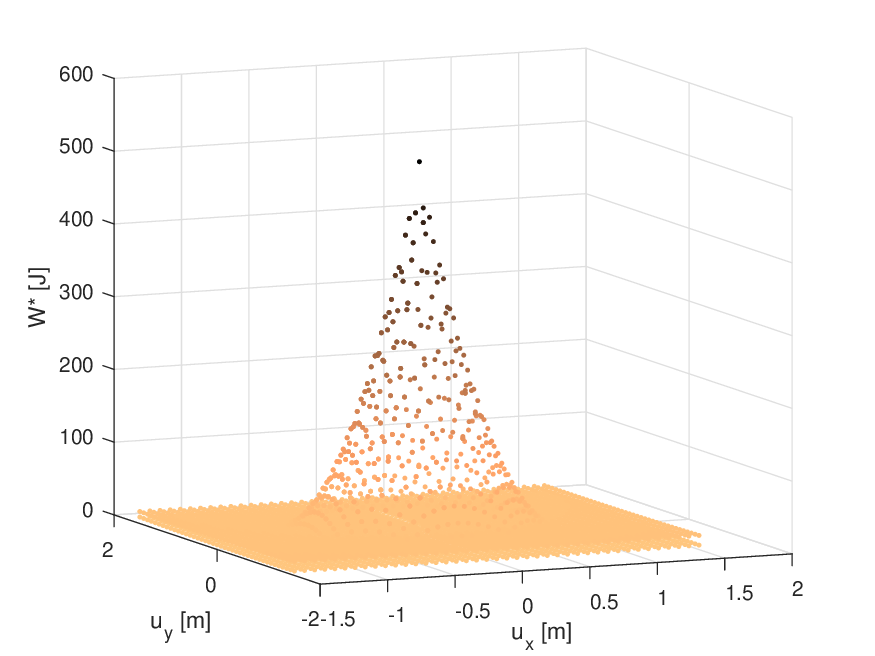}
	\\	b) \\
	\caption{a) Equilibrium configurations $z^*$ (equivalent to the rotation $\alpha$ of the pendulum) and b) equilibrium energies $W^*$  for an array of inputs $(u_x, u_y)$ spanning 150\% of the workspace ($\pm L_0$)}
	\label{fig:Staircase}
\end{figure}
At each point $(u_x, u_y)$ of the grid, multiple solutions $z^*_m$ (a scalar, representing an equilibrium angle $\alpha$ of the pendulum angle), each \HL{with its own    multiplicity index $m$,} were estimated by the numerical solver and are shown in Fig.~\ref{fig:Staircase}-a), where the unstable solutions (the `ceiling' at height $\alpha \approx +\pi/2$) are highlighted in red color. The remaining solutions consist of a `floor' $\alpha \approx -\pi/2$ and a `staircase' which is approximately the graph of \HL{Eq.~}\eqref{EQ:Ex1_atan2}. The `floor' and the `ceiling' in Fig.~\ref{fig:Staircase}-a) are only approximately horizontal due to the presence of a small but non-zero value of $k_{min}$ through which the agent still influences the pendulum at a distance ($d>0$). Figure~\ref{fig:Staircase}-b) shows the energies relative to the same solutions, the higher energies (cone-like shape in the middle) are relative to the equilibria on the `staircase'. It is evident how the amount of elastic energy increases as the agent approaches the origin, due to the compression ($d\ll 0$) of the spring ($k_c(d)\approx k_{max}$).
\HL{
	Note that for very high stiffness, i.e. $k_{max}\to\infty$, the critical point approaches the origin, i.e. $mg/2k_{max}\to 0$.
}

\subsubsection{Optimal Path on Graphs}
In previous section, the the problem of optimally manipulating the pendulum in Fig.~\ref{fig:Ex1} was framed as an optimal path on the equilibrium manifold and solved specifically on the `staircase', one of the possible branches of the equilibrium manifold, analytically  defined by \eqref{EQ:Ex1_atan2}.

In real or numerical scenarios, we rarely have the luxury of an equation such as \eqref{EQ:Ex1_atan2}, locally defining a branch of the equilibrium manifold. We might just have samples of solutions, either from experiments or numerical approximations, and there might be no indication of which branch they belong to (or prior knowledge of the number of branches). In fact, the plots in Fig.~\ref{fig:Staircase} resemble this situation: all we have are \textit{disconnected solutions}, each representing an equilibrium point in and open $\mathcal Z\times \mathcal U\in\mathbb R^{N\times K}$. We shall here only focus on the stable equilibria: numerically, we can always test the Hessian $\det (\nabla_{\V z\V z}W)$ while, experimentally, we will not be able to even witness unstable equilibria (just like we never witness a tossed coin landing and standing on its edge).

To capture the {topology} of the problem, we propose working on \textit{graphs}. 
Referring to textbooks such as \cite{Edelsbrunner Harer 2022} for details, here we simply recall the following basic definition: an undirected graph is a pair $\mathcal G=(\mathcal V, \mathcal E)$ where $\mathcal V$ is a set of vertices and $\mathcal E$ is a set of edges, where each element $e_{ij}\in \mathcal E$ consists of an (unordered) pair of vertices $v_i, v_j\in\mathcal V$. Two vertices are connected (i.e. `neighbors' in some sense, to be specified) if there is an edge connecting them.

\begin{figure}[h!]
	\centering
	\includegraphics[scale = 0.9]{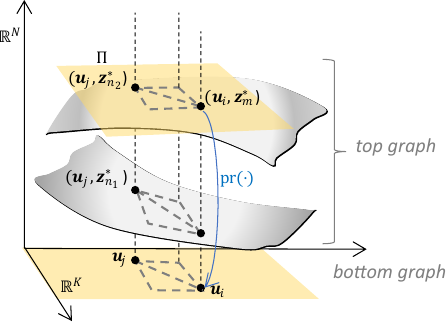}
		\caption{Construction of connectivity in the top-graph' (in full configuration space $\mathbb R^N\times \mathbb R^K$) given the connectivity in the `bottom-graph' (in control space $\mathbb R^K$).
		%The sketch depicts (parts of) two branches of the EM lying `above' the control point $\V u_j\in \mathbb R^K$ and relative to different multiplicities ($\V z_{n_1}\neq \V z_{n_2}$). Not shown in this figure is the fact that these two branches might actually reconnect above some set of \textit{critical} points $\V u_c\in\mathbb R^K$ but, for this to happen, the Hessian $\nabla_{\V z\V z}^2W$ must be rank-deficient when evaluated at $(\V z_c, \V u_c)$ and condition \eqref{EQ:max_rank_Wzz} no longer holds.
	}
	\label{fig:Graph_Bundle}
\end{figure}

With reference to Fig.~\ref{fig:Graph_Bundle}, given the control nature of the problem, we assume that an initial graph, referred to as `bottom graph' $\mathcal G_{bot}$, is available. In practice, this corresponds to having an array (often a regular grid) of $S$ control inputs $\{\V u_s \in \mathcal U\}_{s=1}^S$. These sampling locations constitute the set of vertices for the given bottom graph $\mathcal G_{bot}$. With reference to the pendulum problem in Fig.~\ref{fig:Ex1}, the inputs $\V u_s$ are taken to be a regular grid of a region of the workspace and two vertices are considered connected if they are neighbors in this regular grid. 

\HL{The sketch in Figure~\ref{fig:Graph_Bundle} depicts (parts of) two branches of the EM lying `above' the control point $\V u_j\in \mathbb R^K$ and relative to different multiplicities ($\V z_{n_1}\neq \V z_{n_2}$). Not shown in this figure is the fact that these two branches might actually reconnect above some set of \textit{critical} points $\V u_c\in\mathbb R^K$ but, for this to happen, the Hessian $\nabla_{\V z\V z}^2W$ must be rank-deficient when evaluated at $(\V z_c, \V u_c)$ and condition \eqref{EQ:max_rank_Wzz} no longer holds.}
In particular,  Figure~\ref{fig:Graph_Bundle} shows two vertices $\V u_i, \V u_j\in\mathcal U$ of the bottom graph and an edge (diagonal, dashed line) connecting them. The other edges (forming parallelogram) represent connections to other neighbors and this connectivity pattern is used to tile up the whole grid.

The objective is to construct a `top graph' $\mathcal G_{top}$ given the bottom one $\mathcal G_{bot}$. More specifically, given isolated solutions which constitute the vertices of $\mathcal G_{top}$ and which naturally project to vertices on $\mathcal G_{bot}$, the objective is to re-construct the connectivity on $\mathcal G_{top}$ based on the connectivity on $\mathcal G_{bot}$. 

To each vertex $\V u_i$ on $\mathcal G_{bot}$ there correspond possibly multiple vertices, denoted as $(\V z^*_m, \V u_i)$, in $\mathcal G_{top}$, where $m=1, 2, \dots, $ is the multiplicity of solutions of \eqref{EQ:equil}. Since each solution $\V z^*_m$ is found only after the control input $\V u_i$ is specified,  there exists a \textit{natural projection} 
\begin{equation}\label{EQ:projection}\nonumber
	\text{pr}: (\V z^*_m, \V u_i) \mapsto \V u_i   .
\end{equation}

Given a solution (vertex) $(\V z^*_m, \V u_i)$ on the top graph and control $\V u_j$ connected to $\V u_i$ on the bottom graph, we would like to establish which of the solutions $(\V z^*_n, \V u_j)$ projecting to $\V u_j$ should be connected to $(\V z^*_m, \V u_i)$. In other words, which of the solutions $(\V z^*_n, \V u_j)$ lies on the same branch as $(\V z^*_m, \V u_i)$.

Assuming a relatively dense grid, the idea is to numerically approximate a branch with a tangent space at $(\V z^*_m, \V u_i)$. Therefore, of all possible solutions $\V z^*_n$ projecting down on the same $\V u_j$, we should consider the one for which $\V z^*_n$ is `closest' to \HL{$\V z^*_m +\delta \V z_{j, i,m}$} where, from \HL{Eq.}\eqref{EQ:1st_order}, 
\begin{equation}\nonumber
	\delta \V z_{j, i,m} = -  (\nabla_{\V z\V z}^2W|_{*mi})^{-1}\nabla^2_{\V u\V z}W|_{*mi}(\V u_j-\V u_i),
\end{equation}
where the notation $|_{*mi}$ denotes that the gradients are evaluated at $(\V z_m^*, \V u_i)$.

This is sketched in Fig.~\ref{fig:Graph_Bundle}), where the plane $\Pi$ is tangent at $(\V z^*_n, \V u_j)$ to the $m$-th branch, according to \HL{Eq.}\eqref{EQ:1st_order}. 

\REMARK{
[{\bf Fiber-wise Distance}]
The reasoning above made use of the concept of `closeness' between $\V z^*_n$ and $\V z^*_m +\delta \V z$. This distance is evaluated on the {`fiber' above $\V u_j$} and this distance is inherent to the mechanical configuration of the space $\mathcal Z\HL{\subset}\mathbb R^N$. We assume a fiber-wise distance, perhaps multiple ones, can always for a given problem at hand. In the specific case of the pendulum, it simply refers to difference in orientation $\alpha$ between two configurations and such a {difference is evaluated \textit{modulo} $2\pi$}.
}

Once a connectivity is built on $\mathcal G_{top}$, each edge can be assigned  a \textit{non-negative weight}: if two vertices  $(\V z^*_n, \V u_j)$ and $(\V z^*_m, \V u_i)$ of the top graph $\mathcal G_{top}$ are connected, a non-negative weight $w_{j,i,m}$ is defined as the squared-Hessian:
\begin{equation}\nonumber
		w_{j,i,m} := (\V u_j-\V u_i)^T\V G|_{*mi}^2(\V u_j-\V u_i)   ,
\end{equation}
where $\V G|_{*mi}$ is the Hessian defined in \eqref{EQ:Gm(u)} and evaluated at  $(\V z_m^*, \V u_i)$.

\begin{figure}[h!]
	\centering
	\includegraphics[width=0.70\linewidth ]{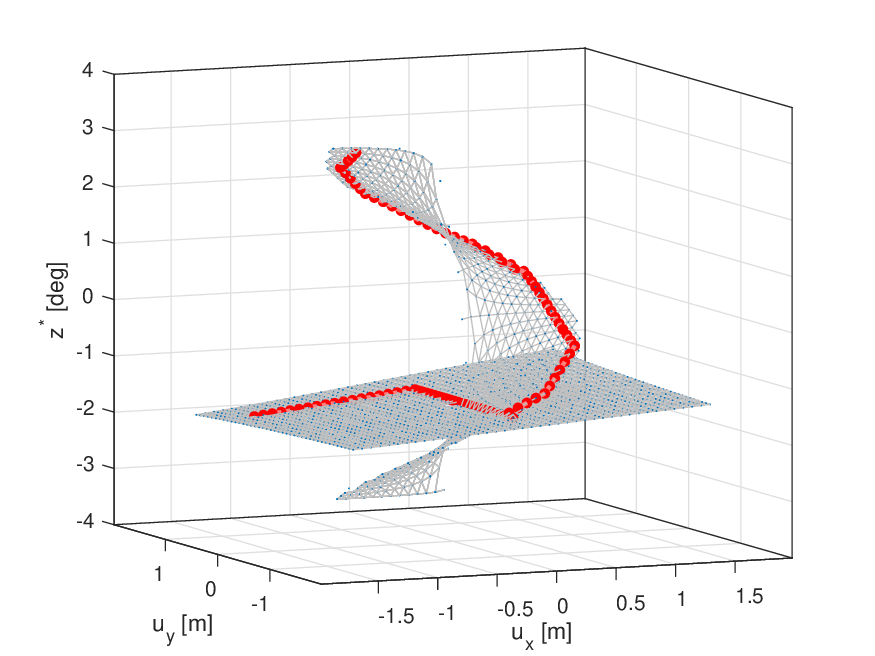}
	\caption{Minimum path over the top-graph $\mathcal G_{top}$ connecting two (arbitrary) points.}
	\label{fig:OptimalPath_on_Graph}	
\end{figure}

Standard search algorithms are available for finding optimal paths on graphs given non-negative weight connectivity. Fig.~\ref{fig:OptimalPath_on_Graph} shows a graph $\mathcal G_{top}$ obtained from numerical solutions following the procedure described above. An optimal path is highlighted connecting two points, chosen to be on different branches.
One noteworthy aspect is that the the only way to go up on the `staircase' (i.e. reversing the pendulum) is to pass through the intersections between `floor' and `staircase'. In this sense, the algorithm is able to navigate a multivalued map, finding passageways in-between branches.

We believe that this might lead to a topological description of the manipulation tasks. 

\HL{
\section{Conclusion}
In this work, we addressed the problem of quasi-static  mechanical manipulation using classical methods from geometric mechanics.
This is done in a two-part paper, as an attempt to bridge geometric insight with data analysis. 
In the first part, more theoretical, classical concepts are reviewed and main assumptions are put in place to describe manipulation as a quasi-static, smooth optimal control problem. 
The second part, more applied, based on a simple but representative problem of mechanical manipulation, illustrates how discrete constructions over the set of numerical solutions derived from the original smooth problem, could be used to provide approximate solutions. 

In the first part, we first recognize how a network of elastically interconnected rigid bodies can be described as a Lagrangian submanifold of the cotangent bundle of its configuration space. 
We also show how the control aspect of the problem can be accounted for by splitting the configuration variables into internal states and input controls. Then, from a second order analysis, we derive a metric directly from the Hessian of the very scalar potential defining the elastic interconnection of the system.

In the second part of the paper, using an elastically- driven, inverted pendulum as a toy model, we show how to regularize mechanical contact to obtain a smooth, albeit highly nonlinear, potential.
Standard numerical approaches can then be used to derive, for example, a discrete set of solutions. Graph theoretical methods are used to construct a weighted graph, using connectivity and metric arguments from the first part of the paper.
In summary, the original manipulation problem reduces to finding a minimum paths over a weighted graph. 

One of the limitations of this study is that it does \textit{not} consider noise, which is inherent not only to experimental datasets  but also to numerical ones, as approximation errors. This limitation hints at the necessity of going beyond equilibrium, in particular a next step will be a near-equilibrium analysis which will allow to consider noise as a form of fluctuation around equilibria.
}

\section*{Acknowledgment}

\HL{The authors would like to thank anonymous reviewers for their insightful and constructive suggestions, which included  references to Tulczyjew's triple.}
\HL{The contribution of F. Cardin to this paper has been realized within the sphere of activities of the GNFM of the INDAM.}
This research is supported by the National Research Foundation, Singapore, under the NRF Medium Sized Centre scheme (CARTIN). Any opinions, findings and conclusions or recommendations expressed in this material are those of the authors and do not reflect the views of National Research Foundation, Singapore. 

\newpage
\appendix\section{Appendix - Tensorial properties of the Hessian}\label{sec:Appendix}

\subsection{Standard formulas from matrix calculus}
\label{sec:Appendix_Hessian}
Consider the following (scalar) potential 
\begin{equation}\nonumber
	f:\R^n\to \R.
\end{equation}
The Hessian at a point $\V x\in\mathbb R^n$ is defined as 
\begin{equation}\nonumber
	\nabla_{\V x\V x}^2 f := \nabla_{\V x}\nabla_{\V x}^T f  .
\end{equation}
Consider a change of coordinates 
\begin{equation}\nonumber
	\V x=\V X(\Vtilde x) 
	\hspace{1cm} \HL{\Longrightarrow} \hspace{1cm} 
	\delta \V x = \V J \delta \Vtilde x,
\end{equation}
where $\V X:\mathbb R^n\to \mathbb R^n$ is locally smooth and invertible, and the new scalar function and $\V J$ its Jacobian. 

In new coordinates system, the potential is defined as 
\begin{equation}\nonumber
	\tilde f(\Vtilde x):=f(\V X(\Vtilde x))
	\hspace{1cm} \HL{\Longrightarrow} \hspace{1cm}
	d\tilde f=df \V J,
\end{equation}
where differentials $df, d\tilde f$ are intended as row vectors, it follows that \HL{in column vector notation, via the \textit{nabla }operator, }  forces transform as
\begin{equation}\label{key}\nonumber
	\nabla_{\Vtilde x}\tilde f = \V J^T\nabla_{\V x}f
\end{equation}
The Hessian in new coordinates is given as 
\begin{equation}\label{eq:Hessian_change_of_coord}
	\nabla^2_{\Vtilde x\Vtilde x}\tilde f = \V J^T \nabla_{\V x\V x}f \V J 
	+\underbrace{
		\MAT{\frac{\partial \V J}{\partial \HL{\tilde x_{1}}}\nabla_{\V x}f &\cdots &\frac{\partial \V J}{\partial \HL{\tilde x_{n}}}\nabla_{\V x}f }
	}_{\text{non-covariant}}.
\end{equation}

\REMARK{
The Hessian is covariant if the second term in \HL{ \eqref{eq:Hessian_change_of_coord} }is null.	This is guaranteed at critical points ($\nabla_{\V x}f=0$).
}

Let $\V K:=\nabla^2_{\V x\V x}f$ and $\Vtilde K:=\nabla^2_{\Vtilde x\Vtilde x}\tilde f$ represent stiffness matrices (i.e. Hessians of some potential) in two different coordinate frames. At equilibrium, the following transformation holds
\begin{equation}\label{eq:Hessian_covariance}
	\Vtilde K= \V J^T\V K \V J .
\end{equation}
Suppose $\delta \V x$ is an eigenvector of $\V K$ with eigenvalue $\lambda$, i.e.
\begin{equation}\label{}\nonumber
	\V K \delta\V x = \lambda \delta\V x .
\end{equation}
The corresponding vector $\delta\Vtilde x$ in new coordinates, i.e. $\delta\V x=\V J\delta\Vtilde x$, is an eigenvector (with same eigenvalue $\lambda$) only if the following {\em orthonormality condition} holds
\begin{equation}\label{eq:ortho}
	\V J^{-1}=\V J^T  .
\end{equation}

{\small
	\textbf{Proof:}
	\begin{equation}\label{}\nonumber
		\begin{matrix}
			\V K \delta \V x 		& =&  \lambda \delta \V x, \\
			\V K \V J\delta \Vtilde x &=& \lambda \V J \delta \Vtilde x, \\
			{\V J^{-1}} \V K \V J\delta \Vtilde x &=& \lambda  \delta \Vtilde x, \\
			\Vtilde K \delta \Vtilde x &=& \lambda  \delta \Vtilde x .\\
		\end{matrix}
	\end{equation}
}

\REMARK{
If Hessians transform covariantly \eqref{eq:Hessian_covariance} and the orthonormality condition \eqref{eq:ortho} holds then the squared-Hessian $\V K^2$ (in fact, any power) will transform covariantly	
	\begin{equation}\label{eq:squared_Hessian_covariance}
		\Vtilde K^2= (\V J^T\V K \V J)^2= \V J^T\V K^2 \V J.
	\end{equation}
}
\subsubsection*{Invariance for our problem}
To fit our problem, consider the following renaming of variables
$$\V x \equiv[z_1, \dots, z_N, u_1, \dots, u_K]^T\equiv\MAT{\V z \\ \V u},$$
$$f(\V x)\equiv W(\V z, \V u) .$$
Note that 
$$\nabla_{\V x}f=\MAT{\nabla_{\V z}W \\ \nabla_{\V u}W}=\MAT{\V 0 \\ \nabla_{\V u}W}.$$
The following {control-affine change of coordinates} renders the non-covariant term in \eqref{eq:Hessian_change_of_coord} null
\begin{equation}\nonumber
	\V X(\Vtilde x)\equiv \MAT{\V z\\ \V u}=\MAT{\V \phi(\Vtilde z)+\V B(\Vtilde z)\Vtilde u \\ \V A(\Vtilde z) \Vtilde u}.
\end{equation}
Note that this leads to a Jacobian  which is independent of controls $\Vtilde u$, i.e.
$$\V J=\MAT{\nabla_{\Vtilde z}\V \phi(\Vtilde z)+ \nabla_{\Vtilde z}\V B(\Vtilde z) \\ \nabla_{\Vtilde z}\V A(\Vtilde z)}
\Longrightarrow\frac{\partial \V J}{\partial \tilde u_k}=\V 0 .
$$

\subsection{Squared-Hessian as a natural metric on the Lagrangian submanifold}
\label{sec:Appendix_squaredHessian}
In this section, we show how to construct an `intrinsic' metric on the Lagrangian submanifold $\Lambda_m$.
For the convenience of the reader, we recall that $W{_m^*}(\V u)=W(\V z_m(\V u),\V u)$ and remark that $\Lambda_m$ is identificable ($\simeq$) by isomorphism with $\mathbb R^K$, then $\iota\simeq dW{_m^*}$ and $T\iota\simeq T dW{_m^*}$. 
%We denote also by $\simeq$ the standard involution between $T^*T\mathbb R^K$ and $TT^*\mathbb R^K$.
Here,  we shall also use $\V f(\V u){=-dW_m^*}$ as a shorthand for $\V f_{ctrl}(\V z_m(\V u), \V u)$ defined in \eqref{EQ:f_eq}, {where $\V z_m(\V u)$ is a well defined function} whenever the rank condition \eqref{EQ:max_rank_Wzz} is satisfied. 

$$ \begin{tikzcd}
	\ARRAY{\mathbb R^K\\ \V u } 
		\arrow[r, shift left=1.5ex, "\simeq" description, phantom]
		\arrow[r, mapsto, shift right=1.5ex, "\text{id}"]
	& \ARRAY{\Lambda_m \\ \V u  } 
		\arrow[r, shift left=1.5ex, "\iota", hookrightarrow]
		\arrow[r, mapsto, shift right=1.5ex]
	& \ARRAY{T^*\mathbb R^K\\   (\V u, \V f(\V u))} 
		\arrow[r, shift left=1.5ex, "\pi_{\mathbb R^K}"]
		\arrow[r, mapsto, shift right=1.5ex]
	& \ARRAY{\mathbb R^K \\ \V u .} 
\end{tikzcd} $$

For $TT^*\mathbb R^K$ there exist the following {\em tangent} {($T_\pi{_{\mathbb R^K}}$)} and { {\em vertical} ($V_\pi{_{\mathbb R^K}}$)} fibrations:

$$
\begin{tikzcd}
	 &  & \ARRAY{T\mathbb R^K \\ (\V u, \delta \V u)}  \\
	\ARRAY{T\Lambda_m \\ (\V u, \delta \V u)} 
	\arrow[r, shift left=1.5ex, "T\iota"]
	\arrow[r, mapsto, shift right=1.5ex]
	& 
	\ARRAY{TT^*\mathbb R^K \\ (\V u, \V f(\V u); \delta \V u, \V f'(\V u)\delta \V u)} 
	\arrow[ur, "T_{\pi_{\mathbb R^K}}"]
	\arrow[dr, "V_{\pi_{\mathbb R^K}}"']
	\\
	 &  & \ARRAY{T^*\mathbb R^K \\ (\V u, \V f'(\V u)\delta \V u).} \\
\end{tikzcd}
$$

\HL{
In order to gain an intrinsic meaning of the above framework, we insert our setting in the Tulczyjew's triple:}

$$
\begin{tikzcd}
    \ARRAY{T^*T^*M \\ (\V q, \V p, -\Vdot p, \Vdot q)}     
        & \ARRAY{TT^*M \\ (\V q, \V p, \Vdot q, \Vdot p)} 
        \arrow[l, shift right=1.5ex, "\beta_M"']
        \arrow[l, mapsto, shift right=-1.5ex]
        \arrow[r, shift right=-1.5ex, "\alpha_M"]
        \arrow[r, mapsto, shift right=1.5ex]
        &    \ARRAY{T^*TM \\ (\V q, \Vdot q, \Vdot p, \V p)} 
\end{tikzcd}
$$

\HL{It was introduced to clarify globally the Legendre transformation; for recent literature on it, see e.g. \cite{Zajac and Grabowska 2016} and its rich bibliography. Operatively, in our case:
}
$$
\begin{tikzcd}%[row sep=large]
\ARRAY{T^*T\R^K \\ (\V u,\delta\V u, \V f(\V u), \V f'(\V u)\delta\V u)} 
\arrow[r, shift right=-1.5ex, "\alpha^{-1}_{\R^K}"]
\arrow[d, "\pi_{T\R^K}"']
&
\ARRAY{TT^* \R^K \\ (\V u,\V f'(\V u)\delta\V u, \delta\V u, \V f(\V u) )} 
\arrow[d, "\tau_{T^*\R^K}"]
\\
\ARRAY{T\R^K \\ (\V u,\delta\V u)}
&
\ARRAY{T^*\R^K \\ (\V u,\V f'(\V u)\delta\V u)}
\end{tikzcd}
$$

On the tangent and cotangent bundles, the Euclidean metric is defined, respectively, as 
$$\begin{tikzcd}
	\ARRAY{	{\mathbb  I}: & T\mathbb R^K \times_{\mathbb R^K} T\mathbb R^K \\ 
		& (({\V u}, \delta \V u_1), ({\V u}, \delta \V u_2))}
		\arrow[r, shift left=1.5ex, ]
		\arrow[r, mapsto, shift right=1.5ex]
	& \ARRAY{\mathbb R \hfill \\ {\mathbb  I}(\delta \V u_1, \delta \V u_2):=\delta_{ij}\delta  u_1^i\delta  u_2^j, }
\end{tikzcd}$$ 			%
$$\begin{tikzcd}
	\ARRAY{	{\mathbb  I}^{-1}: & T^*\mathbb R^K \times_{\mathbb R^K} T^*\mathbb R^K \\ 
		& (({\V u}, \delta \V \varphi_1), ({\V u}, \delta \V \varphi_2))}
	\arrow[r, shift left=1.5ex, ]
	\arrow[r, mapsto, shift right=1.5ex]
	& \ARRAY{\mathbb R \hfill \\ {\mathbb  I}^{-1}(\delta \V \varphi_1, \delta \V \varphi_2):=\delta^{ij}\delta \varphi_{1i}\delta \varphi_{2j}.}
\end{tikzcd}$$			%
An arbitrary curve on $\mathbb R \ni s \mapsto \V u(s) \in \mathbb R^K\simeq_{\text{id}}\Lambda_m$, together with its generated tangent vectors, can be lifted on $T\Lambda_m$ as follows

\begin{equation}\label{eq:Vertical_pullback}
\begin{tikzcd}
	\ARRAY{T\Lambda_m \\ (\V u, \Vdot u)} 
	\arrow[r, shift left=1.5ex, "T\iota"]
	\arrow[r, mapsto, shift right=1.5ex]
	& 
	\ARRAY{TT^*\mathbb R^K \\ (\V u, \V f(\V u);  \Vdot u, \V f'(\V u)  \Vdot u)} 
	\arrow[r, shift left=1.5ex, "V_{\pi_{\mathbb R^K}}"]
	\arrow[r, shift right=1.5ex, mapsto]
	& \ARRAY{T^*\mathbb R^K \\ (\V u, \V f'(\V u)\Vdot u).}  
\end{tikzcd}
\end{equation}
%which \HL{ allows  defining the metric $\mathfrak g_{\Lambda_m}$} as follows

{By the following \textit{vertical} pull-back from $T^*\mathbb R^K$ to $T\Lambda_m$ of the Euclidean metric $\mathbb I^{-1}$ on $T^*\mathbb R^K$, we are able to define the metric $\mathfrak g_{\Lambda_m}$, finally involving $\V G_m^2$ introduced in \eqref{EQ:Gm(u)}:}

$$\ARRAY{
	\mathfrak g_{\Lambda_m}(\V u) (\Vdot u, \Vdot u)
	& :=& {\mathbb  I}^{-1}(V_{\pi_{\mathbb R^K}}\circ T\iota (\Vdot u ), V_{\pi_{\mathbb R^K}}\circ T\iota (\Vdot u ))\\
%	& {= (V_{\pi_{\mathbb R^K}}\circ T\iota (\Vdot u ))^* \HL{\mathbb  I}^{-1}(\Vdot u, \Vdot u)}\hfill \\
	& =& {\mathbb  I}^{-1} (\V f'(\V u)\Vdot u, \V f'(\V u)\Vdot u) \hfill \\
	& =& {\mathbb  I}^{-1} (\V G_m(\V u)\Vdot u, \V G_m(\V u)\Vdot u) \hfill \\
	& =& \Vdot u^T \V G_m^2 \Vdot u. \hfill
	}
$$

\end{document}